\documentstyle[11pt]{article}

\input amssym.def
\input amssym

\let\got=\frak

\def\filename{factor.tex}

\def\hybrid{\topmargin 0pt      \oddsidemargin 0pt
        \headheight 0pt \headsep 0pt
        \textwidth 17.5cm
        \textheight 25cm
        \voffset=-0.2cm
        \hoffset=-0.4cm
        \marginparwidth 0.0in
        \parskip 5pt plus 1pt   \jot = 1.5ex}
\catcode`\@=11
\def\marginnote#1{}

\newcount\hour
\newcount\minute
\newtoks\amorpm
\hour=\time\divide\hour by60
\minute=\time{\multiply\hour by60 \global\advance\minute by-\hour}
\edef\standardtime{{\ifnum\hour<12 \global\amorpm={am}%
        \else\global\amorpm={pm}\advance\hour by-12 \fi
        \ifnum\hour=0 \hour=12 \fi
        \number\hour:\ifnum\minute<10 0\fi\number\minute\the\amorpm}}
\edef\militarytime{\number\hour:\ifnum\minute<10 0\fi\number\minute}

\def\draftlabel#1{{\@bsphack\if@filesw {\let\thepage\relax
   \xdef\@gtempa{\write\@auxout{\string
      \newlabel{#1}{{\@currentlabel}{\thepage}}}}}\@gtempa
   \if@nobreak \ifvmode\nobreak\fi\fi\fi\@esphack}
        \gdef\@eqnlabel{#1}}
\def\@eqnlabel{}
\def\@vacuum{}
\def\draftmarginnote#1{\marginpar{\raggedright\scriptsize\tt#1}}

\def\draft{\oddsidemargin -0.1truein
        \def\@oddfoot{\sl preliminary draft\ {\tt\filename} \hfil
        \rm\thepage\hfil\sl\today\quad\militarytime}
        \let\@evenfoot\@oddfoot \overfullrule 3pt
        \let\label=\draftlabel
        \let\marginnote=\draftmarginnote
   \def\@eqnnum{{\rm (\theequation)}\rlap{\kern\marginparsep\tt\@eqnlabel}%
\global\let\@eqnlabel\@vacuum}  }

\newcommand{\CC}{{\Bbb C}}

\newcommand{\ZZ}{{\Bbb Z}}
\let\z=\ZZ

\font\fiveeufm=cmmib5  

\font\yao=cmmi12 
\def\bep{{\textfont1=\yao\scriptfont1=\teneuf
\scriptscriptfont1=\fiveeufm
\hbox{$\mathsurround=0pt\displaystyle\epsilon$}}
}

\newdimen\linethick  \linethick=0.4pt
\newdimen\hboxitspace    \hboxitspace=5pt
\newdimen\vboxitspace    \vboxitspace=5pt

\def\fr#1{%
\beq\new
\vcenter{
\hrule height\linethick
           \hbox{\vrule width\linethick
                 \kern\hboxitspace
                 \vbox{\kern\vboxitspace
                       \hbox{$\begin{array}{c}\displaystyle#1
          \end{array}$}%
                       \kern\vboxitspace}%
                 \kern\hboxitspace
                 \vrule width\linethick}%
           \hrule height\linethick}%
\eeq}

\newdimen\Squaresize \Squaresize=14pt
\newdimen\Thickness \Thickness=0.5pt

\def\Square#1{\hbox{\vrule width \Thickness
   \vbox to \Squaresize{\hrule height \Thickness\vss
      \hbox to \Squaresize{\hss#1\hss}
   \vss\hrule height\Thickness}
\unskip\vrule width \Thickness}
\kern-\Thickness}

\def\Vsquare#1{\vbox{\Square{$#1$}}\kern-\Thickness}

\def\numberbysection{\@addtoreset{equation}{section}
        \def\theequation{\thesection.\arabic{equation}}}
\numberbysection

\renewcommand{\theequation}{\thesection.\arabic{equation}}
\newcommand{\l@qq}[2]{\addvspace{2em}
 \hbox to\textwidth{\hspace{1em}\bf #1 \dotfill #2}}

\def\appname{Appendix}
\newcounter{app}
\def\theapp{\Alph{app}}
\def\app{\par
   \addvspace{4ex}
   \@afterindentfalse
  \secdef\@app\@dapp}
\def\@app[#1]#2{\ifnum \c@secnumdepth >\m@ne
        \refstepcounter{app}
        \addcontentsline{toc}{app}{\theapp
        \hspace{1em}#1}\else
      \addcontentsline{toc}{app}{ #1}\fi
   {\parindent \z@ \raggedright
    \Large \bf \appname~\theapp .
   \Large  \bf \hspace{1em}    #2}\nobreak
   \vskip -2ex   \noindent
\setcounter{equation}{0}
\def\theequation{\Alph{app}.\arabic{equation}}}
\def\@dapp#1{%
{\parindent \z@ \raggedright  \bf #1}\par\nobreak}
\def\l@app#1#2{\addpenalty{\@secpenalty}%
   \addvspace{1em plus\p@}%
   \begingroup
   \@tempdima 3em
     \parindent \z@ \rightskip \@pnumwidth
     \parfillskip -\@pnumwidth
     { \bf
     \leavevmode
     #1\hfil \hbox to\@pnumwidth{\hss #2}}\par
     \nobreak
   \endgroup}
\newcounter{sapp}[app]
\def\thesapp{\Alph{app}.\arabic{sapp}}
\def\sapp{\par
   \addvspace{4ex}
   \@afterindentfalse
  \secdef\@sapp\@dsapp}
\def\@sapp[#1]#2{\ifnum \c@secnumdepth >\m@ne
        \refstepcounter{sapp}
        \addcontentsline{toc}{sapp}{\thesapp
        \hspace{1em}#1}\else
      \addcontentsline{toc}{sapp}{ #1}\fi
   {\parindent \z@ \raggedright
    \large \bf \thesapp
   \large  \bf \hspace{1em}    #2}\nobreak
   \vskip 4ex   \noindent
\def\theequation{\Alph{app}.\arabic{equation}}}
\def\@dsapp#1{%
{\parindent \z@ \raggedright  \bf #1}\par\nobreak}
\def\l@sapp#1#2{\addpenalty{\@secpenalty}%
   \begingroup
   \@tempdima 3em
     \parindent \z@ \rightskip \@pnumwidth
     \parfillskip -\@pnumwidth
     { \hspace{1em}
     \leavevmode
     #1 \hfil \dotfill \hbox to\@pnumwidth{\hss #2}}\par \nobreak
     \endgroup}

\def\titlepage{\@restonecolfalse\if@twocolumn\@restonecoltrue\onecolumn
     \else \newpage \fi \thispagestyle{empty}\c@page\z@
        \def\thefootnote{\fnsymbol{footnote}} }

\def\endtitlepage{\if@restonecol\twocolumn \else  \fi
        \def\thefootnote{\arabic{footnote}}
        \setcounter{footnote}{0}}  
\relax

\hybrid

\newtoks\@stequation

\def\subequations{\refstepcounter{equation}%
  \edef\@savedequation{\the\c@equation}%
  \@stequation=\expandafter{\theequation}
  \edef\@savedtheequation{\the\@stequation}
  \edef\oldtheequation{\theequation}%
  \setcounter{equation}{0}%
  \def\theequation{\oldtheequation\alph{equation}}}

\def\endsubequations{%
  \setcounter{equation}{\@savedequation}%
  \@stequation=\expandafter{\@savedtheequation}%
  \edef\theequation{\the\@stequation}%
  \global\@ignoretrue}

\makeatletter
\newdimen\normalarrayskip              
\newdimen\minarrayskip                 
\normalarrayskip\baselineskip
\minarrayskip\jot
\newif\ifold             \oldtrue            \def\new{\oldfalse}
\def\arraymode{\ifold\relax\else\displaystyle\fi} 
\def\eqnumphantom{\phantom{(\theequation)}}     
\def\@arrayskip{\ifold\baselineskip\z@\lineskip\z@
     \else
     \baselineskip\minarrayskip\lineskip1\baselineskip\fi}

\def\@arrayclassz{\ifcase \@lastchclass \@acolampacol \or
\@ampacol \or \or \or \@addamp \or
   \@acolampacol \or \@firstampfalse \@acol \fi
\edef\@preamble{\@preamble
  \ifcase \@chnum
     \hfil$\relax\arraymode\@sharp$\hfil
     \or $\relax\arraymode\@sharp$\hfil
     \or \hfil$\relax\arraymode\@sharp$\fi}}

\def\@array[#1]#2{\setbox\@arstrutbox=\hbox{\vrule
     height\arraystretch \ht\strutbox
     depth\arraystretch \dp\strutbox
     width\z@}\@mkpream{#2}\edef\@preamble{\halign \noexpand\@halignto
\bgroup \tabskip\z@ \@arstrut \@preamble \tabskip\z@ \cr}%
\let\@startpbox\@@startpbox \let\@endpbox\@@endpbox
  \if #1t\vtop \else \if#1b\vbox \else \vcenter \fi\fi
  \bgroup \let\par\relax
  \let\@sharp##\let\protect\relax
  \@arrayskip\@preamble}
\def\eqnarray{\stepcounter{equation}%
              \let\@currentlabel=\theequation
              \global\@eqnswtrue
              \global\@eqcnt\z@
              \tabskip\@centering                      
              \let\\=\@eqncr
              $$%
            \halign to \displaywidth  \bgroup
             \eqnumphantom \@eqnsel
      \hskip\@centering                               
    $\displaystyle  \tabskip\z@ {##}$%
    &\global\@eqcnt\@ne \hskip 2\arraycolsep
         $ \displaystyle  \arraymode{##}$\hfil
    &\global\@eqcnt\tw@ \hskip 2\arraycolsep
         $\displaystyle\tabskip\z@{##}$\hfil
         \tabskip\@centering
    &{##}\tabskip\z@\cr}
\makeatother
\newtheorem{th}{Theorem}
\newtheorem{prop}{Proposition}[section]           
\newtheorem{cor}{Corollary}
\newtheorem{lem}{Lemma}[section]

\def\bea{\begin{eqnarray}}
\def\eea{\end{eqnarray}}

\def\beq{\begin{equation}}
\def\eeq{\end{equation}}
\def\be{\beq\new\begin{array}{c}}
\def\ee{\end{array}\eeq}
\def\bse{\begin{subequations}}                
\def\ese{\end{subequations}}                 %
\def\vbn{$\raisebox{0.1pt}{$\,\hspace{2pt}\llcorner\!$}\leaders
\hrule\hfill\leaders\hrule\hfill\raisebox{0.1pt}{$\!\lrcorner$}\;\,$}
\def\con#1{\mathsurround=0pt
\mathop{\vtop{\ialign{##\crcr$\hfil\displaystyle{#1}\hfil$\crcr
\noalign{\kern1pt\nointerlineskip}\vbn\crcr\noalign{\kern1pt}}}}\limits}
\def\stackreb#1#2{\mathrel{\mathop{#2}\limits_{#1}}}

\def\tr{\triangleright}                                  
\def\tl{\triangleleft}

\def\jo{\mathrel{\mkern-4mu}}
\def\sem{\mathsurround=0pt
\mathrel{\raise1.4pt\hbox{$\scriptscriptstyle>$}}\jo\mathrel\tl}
\def\mes{\mathsurround=0pt
{\mathrel\tr\jo\mathrel{\raise1.4pt\hbox{$\scriptscriptstyle <$}}}}

\def\]{]\raise-2pt\hbox{$_\ast$}}

\def\op#1{\raise-6pt\hbox{$\stackrel{\displaystyle\oplus }{\scriptstyle
#1}$}\;}
\def\oop#1#2{\raise-6pt\hbox{$\stackrel{#2}{\stackrel{\displaystyle\oplus
}{\scriptstyle #1}}$\;}}

\def\al{\alpha}
\def\a{\alpha}
\def\b{\beta}
\def\la{\lambda}

\def\xp{e}
\def\xm{f}
\def\co{\mathsurround=0pt\e\raise-4pt
\hbox{${\!{_{_\cf}\!{\scriptscriptstyle\cal H}}}$}}
\def\dco{\mathsurround=0pt\bep\raise-4pt
\hbox{${\!\!{_{_\cf}\!\!{\scriptscriptstyle\cal H^{\!\ast}}}}$}}
\def\<{\langle}
\def\>{\rangle}
\def\ov{\overline}

\def\hslt{\widehat{\got{sl}}_2}
\def\ca{{\cal A}}
\def\cb{{\cal B}}
\def\ck{{\cal K}}
\def\cu{{\cal U}}
\def\kr{{\cal R}}
\def\cd{{\cal D}}

\def\cf{{\cal F}}
\def\1{1\!\!1}
\def\stack#1#2{\raise0.7pt\hbox{$\mathrel{\mathop{#2}\limits^{#1}}$}}
\def\m{\raise5pt\hbox{$\scriptstyle -1$}}
\def\st#1#2{\raise1pt\hbox{$\stackrel{#1}{#2}$}}
\def\id{\mbox{\rm id}}
\def\fd{\cu\!\mathsurround=0pt
\raise-7pt\hbox{$\stackrel{\displaystyle\otimes}
{\scriptscriptstyle \cal R^{\!\scriptscriptstyle(\!-\!)}}$}\cu}

\def\r#1{(\ref{#1})}
\def\rav#1{\stackrel{\mbox{\tiny {\rm \r{#1}}}}{=}}
\def\ravv#1#2{\stackrel{\mbox{\tiny {\rm \r{#1},\r{#2}}}}{=}}

\def\ot{\otimes}
\def\a{\alpha}    \def\la{\lambda}
\def\b{\beta}  

\def\fp{f_+} \def\fm{f_-}
\def\fpm{f_\pm}
\def\epm{e_\pm}
\def\ep{e_+} \def\em{e_-}

\newcommand{\Uqg}{{U_{q}^{}(\widehat{\frak{g}})}}

\newcommand{\Uqdva}{{U_{q}^{}(\widehat{{sl}}_{2})}}
\newcommand{\sldva}{{\widehat{{sl}}_{2}}}

\let\bn=\be
\let\ed=\ee

\def\sk#1{\left(#1\right)}
\def\skk#1{\bigl(#1\bigr)}

\def\e#1{e^{({#1})}}

\def\T#1{t^{({#1})}}

\def\dz{\underline{dz}}
\def\dzz#1{\frac{dz_{#1}}{z_{#1}}}
\def\dzzz#1{\frac{dz_{#1}}{2\pi iz_{#1}}}
\def\Uep{{U}^+_e}

\def\Ufm{{U}^-_f}

\def\Pep{{P}^+_e}
\def\Pem{{P}^-_e}

\def\Pemp{{P}^\mp_e}
\def\Pfpm{{P}^\pm_f}
\def\Pfp{{P}^+_f}
\def\Pfm{{P}^-_f}
\def\Pepd{{P}^{+*}_e}
\def\Pemd{{P}^{-*}_e}

\def\Pfpmd{{P}^{\pm*}_f}
\def\Pfpd{{P}^{+*}_f}
\def\Pfmd{{P}^{-*}_f}
\def\res#1{\stackreb{#1}{\mbox{\rm res}}}

\def\Pia#1{P^{#1}_{\ca}}
\def\Pib#1{P^{#1}_{\cb}}
\def\tkr{\tilde{\kr}}

\def\dw{\frac{dw}{2\pi iw}}
\def\ant{a}
\def\coun{\varepsilon}
\def\Ev{{\cal E}v}

\begin{document}
\thispagestyle{empty}
\setcounter{page}1

\begin{center}
\bigskip\bigskip
{\Large\bf Factorization of the  Universal ${\cal R}$-matrix
for $U_q(\widehat{sl}_2)$}\\
\bigskip
\bigskip
{\bf
Jintai Ding$^{*}$\footnote{E-mail: Jintai.Ding@math.uc.edu},\
Sergei Khoroshkin$^{\star}$\footnote{E-mail: khor@heron.itep.ru},\
Stanislav
 Pakuliak$^{\star\bullet\circ}$}\footnote{E-mail: pakuliak@thsun1.jinr.ru}\\
\bigskip
$^{*}$
{\it Department of  Mathematical Sciences, University of Cincinnati,\\
 PO Box 210025, Cincinnati, OH 45221-0025,  USA}\\
$^\star$
{\it Institute of Theoretical \& Experimental Physics, 117259 Moscow, Russia}\\
$^\bullet$
{\it Bogoliubov Laboratory of Theoretical Physics, JINR,
141980 Dubna, Moscow region, Russia}\\
$^\circ$
{\it Bogoliubov Institute for Theoretical Physics, 003143 Kiev, Ukraine}
\bigskip
\bigskip
\end{center}

\begin{abstract}
\noindent
The factorization of the universal
${\cal R}$-matrix  corresponding to so called  Drinfeld
Hopf
structure is described on the example of quantum affine algebra
$U_q(\widehat{sl}_2)$.
As a result of factorization procedure we deduce certain differential
 equations on the factors of the universal
${\cal R}$-matrix, which allow to construct uniquely these factors in the
integral form.
\end{abstract}

\footnotesize
\tableofcontents
\newpage
\normalsize

\section{Introduction}

The theory of quantum groups is the
origin
of many group-theoretical
methods for investigation of the quantum integrable models.
The quantum group theory based on the  quantum inverse scattering method
\cite{FT}
was described in  the pioneering works \cite{D1,J}
as the Hopf algebra deformation of the universal enveloping algebras of
contragredient Lie algebras.
 In most  applications the quantum groups
as the Hopf algebras appeared together with $R$-matrices, either in the
form
of numerical matrices or $L$-operators or universal $\kr$-matrices. The
latter are the elements in the completed square of the corresponding
Hopf algebras, which satisfy certain conditions.
A direct consequence of these relations is the fact that
the universal $\kr$-matrix satisfy Yang-Baxter relation and different
$L$-operators and numerical $R$-matrices can be obtained from the
universal
one by specializing to certain representations of the original algebra.

Recently
it became clear \cite{F,ABRR,JKOS} that the deformed algebras
which are behind the integrability of the elliptic models \cite{ABF,Ba}
and their counterparts from the quantum field theories
\cite{ZZ} can be obtained
using twisting procedure
slightly relaxing the coassociativity axiom of the original Hopf
algebras.
The resulting algebras turn out to be the quasi-Hopf algebras
\cite{D2}. The notion of the
universal $\kr$-matrix survive during the twisting procedure although
the
algebra itself can loose some properties. The precise construction of the
twisting element refers to the solution of certain difference equation.
 As
a consequence, the universal $\kr$-matrix in twisted algebra appears
 in this
approach as infinite product of shifted universal $\kr$-matrices
for original quantum affine algebra, which is unobservable for practical
use.

 Another way to get the same result is to use so called 'new realization'
 of quantum affine algebras and their generalizations \cite{D}. This
realization was introduced by Drinfeld in order to show the deformation of
standard loop basis of affine Lie algebras. It
happened to be very useful in representation theory. 'New realization'
possesses its own comultiplication structure (we call it 'Drinfeld
comultiplication'),
 different from standard
comultiplication structure for quantized Kac-Moody algebras. It was
proved in
\cite{KT} that this comultiplication structure can be obtained from the
standard one as a twist by certain factor of universal $\kr$-matrix.
Another
advantage of 'new realizations'
was noted in \cite{JKOS}: their elliptic
analogs can be described with a help of very simple twist.

Enriquez, Felder and Rubtsov \cite{EF,ER} suggested to reverse the
calculations  in \cite{KT}: one can try to describe traditional Hopf
 structure of quantum
affine algebras and their elliptic analogs starting from Drinfeld
comultiplication and its elliptic analog.
It follows from \cite{KT}, that this problem is equivalent to a
Riemann type problem
 of factorization of essential part of the universal $\kr$-matrix for
Drinfeld comultiplication.
They managed to get in this way
an $L$-operator description of elliptic face type algebras and generalized
this approach to the curves of higher genus.

In this paper we develop the ideas of \cite{EF} and solve explicitly the
factorization  problem for quantum affine algebra $\Uqdva$. Our method is
different from \cite{EF} and \cite{ER}. It is based on the use of the
results of \cite{DK,DKP}, where an integral presentation of the universal
$\kr$-matrix for Drinfeld comultiplication was studied. Applying the
projectors which  describe
 the factorization to  this integral presentation
 we
 deduce differential equations for the factors of the universal
$\kr$-matrix
 for $\Uqdva$. These differential
equations have precise unique solution in noncommutative power
 series.
 In particular cases, for instance,
for level one representations, they allow to
describe an evaluation of the universal $\kr$-matrix in
infinite-dimensional
representations. The method is quite general, it can be applied for
the elliptic and Yangian algebras as well,
 though we restrict ourselves to $\Uqdva$ in this paper.

The paper is organized as follows. First, we remind two descriptions of
quantum affine algebra $\Uqdva$ and formulate the main results. Section 3
is devoted to 'new realization' of $\Uqdva$. We describe here its Hopf
structure, the Hopf pairing between two Borel subalgebras in a current
 form and review all known description of the corresponding universal
$\kr$-matrix and
of the pairing tensor. In particular, we remind the results of
\cite{DK,DKP},
where an integral presentation and the differential equation for the
universal $\kr$-matrix was studied.

The main results are proved in Section 4. Here we
first review the projection technique, elaborated in \cite{ER,EF}. Then we
apply it to the differential equation from the previous section and deduce
the differential equation for the factors of traditional $\kr$-matrix.
Here the screening operators naturally appear. In section 5 we compare our
technique with approach from \cite{EF,ER}. We see here that our achievement
can be manifested as a precise calculation of certain multiple integrals,
which can be reduced to a deduction of certain combinatorial identity.
The proof of this identity was found by A.~Okounkov.
 A technique appropriate for direct calculation
of these multiple integrals in this manner, is given
in the
Appendix. It is based on the use of orthogonal bases in Borel subalgebras
and on the intensive use of the pairing technique.

\section{Preliminaries and main results}
\label{preliminary}

\subsection{Two descriptions of the quantum affine algebra $\Uqdva$}

{\bf 1. $\Uqdva$ as quantized Kac-Moody algebra.}\\
Quantum affine algebra $\Uqdva$
is an associative algebra generated by the elements, $e_{\pm\al_i}$,
$k_{\al_i}^{\pm1}$, $d$, $i=0,1$
subjected to the commutation relations:
\be\label{re1}
[d,e_{\pm\al_i}]\,=\,\pm\,\delta_{i0}e_{\pm\al_i},\quad
k_{\al_i}e_{\pm\al_j}k^{-1}_{\al_i}\,=\,q_i^{\pm a_{ij}}e_{\pm\al_j},\quad
[e_{\al_i},e_{-\al_j}]\,=\,
\delta_{ij}\frac{k_{\al_i}-k^{-1}_{\al_i}}{q_i-q_i^{-1}}
\ee
and also
the Serre relations
\be\label{re2}
e^3_{\pm\al_i}e_{\pm\al_j}+[3]_qe^2_{\pm\al_i}e_{\pm\al_j}e_{\pm\al_i}+
[3]_qe_{\pm\al_i}e_{\pm\al_j}e^2_{\pm\al_i}+e_{\pm\al_j}e^3_{\pm\al_i}
=0,\quad
i\neq j,
\ee
where $[n]_q=\frac{q^n-q^{-n}}{q-q^{-1}}$ is Gauss $q$-number and
$a_{ij}=(\al_i,\al_j)$ is Cartan matrix of the affine algebra $\sldva$
$$
a_{ij}=\sk{\begin{array}{cc}2&-2\\-2&2\end{array}}.
$$
The element $c$ given by the relation
\be\label{}
q^c\;\equiv\; k_{\al_0}k_{\al_1}
\ee
is central and its value on the particular representation
is called `level'.

One of the possible
Hopf structures is given by the formulas:
\be\label{copr}
\Delta(e_{\al_i})=e_{\al_i}\otimes 1+k_{\al_i}\otimes e_{\al_i},\quad
\Delta(e_{-\al_i})=1\otimes e_{-\al_i}+e_{-\al_i}\otimes k^{-1}_{\al_i}\\
\Delta(k_{\al_i})\,=\,k_{\al_i}\otimes k_{\al_i}\ ,\quad
\Delta(d)\,=\,d\otimes 1+1\otimes d\\
\coun(e_{\pm\al_i})=0,\quad  \coun(k^{\pm1}_{\al_i})=1,\quad
\coun(d)=0,\\
\ant(e_{\al_i})=-k^{-1}_{\al_i}e_{\al_i},\quad
\ant(e_{-\al_i})=-e_{-\al_i}k_{\al_i},\quad
\ant(k^{\pm1}_{\al_i}) =  k^{\mp1}_{\al_i},\quad
\ant(d)=-d,
\ee
where $\Delta$, $\coun$ and $\ant$ are comultiplication, counit and
antipode maps respectively.
\bigskip

\noindent
{\bf 2. The 'new realization' of $\Uqdva$.}\\ The  so
called new realization of quantum affine algebras was introduced by
V.~Drinfeld
in \cite{D1}.  In this description the algebra $\Uqdva$ is generated by the
infinite set of generators $d$, $q^c$, $k^{\pm1}$, $e_{n}$, $f_{n}$, $n\in\ZZ$,
$a_{\pm m}$, $m\geq0$ subjected to quadratic commutation relations, which are
 given as formal power series identities on their generating functions
\be\label{currents} e(z)\,=\,\sum_{n\in\z}e_{n}z^{-n}\ ,
\quad f(z)\,=\,\sum_{n\in\z}f_{n}z^{-n}\ , \\
\psi^\pm(z)\,=\,\sum_{n\geq0}\psi^\pm_{n}z^{\mp n}
=k^{\pm1}
\,
\exp\sk{\pm(q-q^{-1})\sum_{n>0}a_{\pm n}z^{\mp n}}\ ,
\ee
as follows:
\bse\label{ccc}
\be\label{-1}
[q^c,\mbox{\rm everything}]\,=\,0\ ,
\ee
\be\label{0}
x^d e(z) x^{-d} = e(xz),\quad
x^d f(z) x^{-d} = f(xz),\quad
x^d \psi^\pm(z) x^{-d} = \psi^\pm(xz)\ ,
\ee
\be
(z-q^{2}w)e(z)e(w)=(q^{2}z-w)e(w)e(z)\ ,
\label{1}
\ee
\be
(z-q^{-2}w)f(z)f(w)=(q^{-2}z-w)f(w)f(z)\ ,
\label{2}
\ee
\be
\frac{(q^{c/2}z-q^{2}w)}
{(q^{2+c/2}z-w)}\psi^+(z)e(w)=
e(w)\psi^+(z)\ ,
\label{3a}
\ee
\be
\psi^-(z)e(w)=\frac{(q^{2-c/2}z-w)}
{(q^{-c/2}z-q^{2}w)}
e(w)\psi^-(z)\ ,
\label{4}
\ee
\be
\frac{(q^{-c/2}z-q^{-2}w)}
{(q^{-2-c/2}z-w)}\psi^+(z)f(w)=
f(w)\psi^+(z)\ ,
\label{5}
\ee
\be
\psi^-(z)f(w)=\frac{(q^{-2+c/2}z-w)}
{(q^{c/2}z-q^{-2}w)}
f(w)\psi^-(z)\ ,
\label{6}
\ee
\be
\frac{(z-q^{{2}-c}w)(z-q^{{-2}+c}w)}
{(q^{{2}+c}z-w)(q^{{-2}-c}z-w)}
\psi^+(z)\psi^-(w)=\psi^-(w)\psi^+(z)\ ,
\label{7}
\ee
\be
\psi^\pm(z)\psi^\pm(w)=\psi^\pm(w)\psi^\pm(z)\ ,
\label{7a}
\ee
\be
[e(z),f(w)]=\frac{1}{q-q^{-1}}\left( \delta(z/q^cw)
\psi^+(zq^{-c/2})-\delta(zq^c/w)\psi^-(wq^{-c/2})\right)\ ,
\label{10}
\ee
\ese
where $\delta(z)=\sum_{n\in\ZZ}z^n$.

The coalgebraic structure of the algebra $\Uqdva$ which is equivalent to
\r{copr} cannot be formulated in total currents \r{currents}. Let
\be\label{h-c}
\ep(z)=
\oint \frac{dw}{2\pi w}\ \frac{e(w)}{1-w/z}=
\sum_{k\geq 0}e_{k}z^{-k}\ ,\quad
\em(z)=-
\oint \frac{dw}{2\pi w}\ \frac{e(w)z/w}{1-z/w}=-\sum_{k<0}e_{k}z^{-k}\\
\fp(z)=
\oint \frac{dw}{2\pi w}\ \frac{f(w)w/z}{1-w/z}
=\sum_{k> 0}f_{k}z^{-k}\ ,\quad
\fm(z)=-
\oint \frac{dw}{2\pi w}\ \frac{e(w)}{1-z/w}=-\sum_{k\leq 0}f_{k}z^{-k}
\ee
be the half-currents. Then the coalgebraic structure of
$\Uqdva$ in new realization
have the form (see \cite{KLP}, where the formulas \r{cop} have been
proved in slightly different situation)
\bse\label{cop}
\be\label{copa}
\Delta\sk{e_\pm(z)}=e_\pm(z)\ot1+\sum_ {k\geq0}(-q)^{k}(q-q^{-1})^{2k}
f^k_\pm(zq^2)\ \psi_\pm\sk{zq^{\mp\frac{c_1}{2}}}
\ot e^{k+1}_\pm\sk{zq^{\mp{c_1}}},
\ee
\be\label{copb}
\Delta\sk{f_\pm(z)}=1\ot f_\pm(z)+\sum_ {k\geq0}(-q)^{-k}(q-q^{-1})^{2k}
f^{k+1}_\pm\sk{zq^{\pm{c_2}}}\ot \psi_\pm\sk{zq^{\pm\frac{c_2}{2}}}\
 e^k_\pm(zq^2),
\ee
\be\label{copc}
\Delta\sk{\psi_\pm(z)}=\sum_ {k\geq0}(-1)^{k}[k+1]_q(q-q^{-1})^{2k}
f^{k}_\pm\sk{zq^{2\mp\frac{c_1}{2}\pm\frac{c_2}{2}}}
\psi_\pm\sk{zq^{\pm\frac{c_2}{2}}}
\ot \psi_\pm\sk{zq^{\mp\frac{c_1}{2}}}
 e^k_\pm\sk{zq^{2\mp\frac{c_1}{2}\pm\frac{c_2}{2}}}.
\ee
\ese
where $c_1=c\ot1$ and $c_2=1\ot c$.

\noindent
{\bf 3. The connection between two descriptions.}\\
In the original paper \cite{D1}  the Chevalley generators were presented
in terms of the current generators. In the case under consideration
they have the form:
\be\label{Dr-for}
e_{\al_1}=e_0,\quad e_{-\al_1}=f_0,\quad
e_{\al_0}=q^c f_1 k^{-1},\quad
e_{-\al_0}=k e_{-1} q^{-c},\quad k_{\al_1}=k,\quad  k_{\al_0}=q^ck^{-1}.
\ee
The inverse formulas
appeared lately  in \cite{B,KT1} in general
case of arbitrary quantum affine algebra
 $\Uqg$ (see also \cite{Da} for $\Uqdva$).
It turns out the the generators in `new realization'
coincide up to central elements with   Cartan-Weyl generators of quantum
affine algebra. We describe  shortly this construction in case under
consideration.

First of all, define the generators
$$
e_\delta\;\stackrel{\mbox{\tiny def}}{=}\;
e_{\al_1}e_{\al_0}-q^2e_{\al_0}e_{\al_1}\ ,
\qquad
e_{-\delta}\;\stackrel{\mbox{\tiny def}}{=}\;
 e_{-\al_0}e_{-\al_1}-q^{-2}e_{-\al_1}e_{-\al_0}\ .
$$
They satisfy the commutation relation
$$
[e_\delta,e_{-\delta}]\,=\,\frac{q^c-q^{-c}}{q-q^{-1}}\ .
$$

Now we define the Cartan-Weyl basis of $U_q(\hslt)$.
For $n\geq 0$ we define the generators which correspond to all real roots
$\pm\al_1\pm n\delta$,
$$
E_{\al_1+n\delta}\;\stackrel{\mbox{\tiny def}}{=}\;
[2]_q^{-n}\Big(\mbox{ad}\,e_{\delta}\Big)^ne_{\al_1}\ ,\quad
E_{\al_1-n\delta}\;\stackrel{\mbox{\tiny def}}{=}\;
q^{cn}\,[2]_q^{-n}\Big(\mbox{ad}\,e_{-\delta}\Big)^ne_{\al_1}\ ,
$$
$$
E_{-\al_1+n\delta}\;\stackrel{\mbox{\tiny def}}{=}\;
q^{-cn}(-[2]_q)^{-n}\Big(\mbox{ad}\,e_{\delta}\Big)^ne_{-\al_1}\ ,
\quad
E_{-\al_1-n\delta}\;\stackrel{\mbox{\tiny def}}{=}\;
(-[2]_q)^{-n}\Big(\mbox{ad}\,e_{-\delta}\Big)^ne_{-\al_1}\ .
$$
Define also auxiliary generators
$$
E_{n\delta}\,=\,[E_{\al_1},E_{-\al_1+n\delta}]\ ,\quad
E_{-n\delta}\,=\,[E_{\al_1-n\delta},E_{-\al_1}]\ ,\quad
n>0
$$
related to imaginary roots $\pm n\delta$ generators $a_{\pm n}$
$$
\pm(q-q^{-1})\sum_{n=0}^\infty q^{\pm\frac{cn}{2}} E_{\pm n\delta}
z^{\mp n}=k_{\al_1}^{\pm1}\exp\sk{\pm(q-q^{-1})\sum_{n>0}a_{\pm n}z^{\mp n}}.
$$

The identification with Drinfeld's generators reads as follows
$$
e_n=E_{\al_1+n\delta}\ ,\quad f_n=E_{-\al_1+n\delta}\ ,\quad
\forall\,n\in\ZZ\ ,\quad \psi^\pm_0=k^{\pm1}_{\al_1}\ ,
$$
$$
E_{n\delta}\,\equiv\,\frac{q^{-\frac{cn}{2}}}{q-q^{-1}}\,\psi^+_n
\ ,\quad n>0,\quad
E_{n\delta}\,\equiv\,-\,\frac{q^{-\frac{cn}{2}}}{q-q^{-1}}\,\psi^-_n
\ ,\quad n<0\ .
$$

\subsection{The main results}

Let us remind that a universal $\kr$-matrix for the quantum affine
algebra $\Uqg$ is an element in some completion of
$\Uqg\ot\Uqg$ which  satisfies
\bse\label{u-con}
\be\label{u-con1}
\kr \Delta(x)=\Delta'(x) \kr,\quad \forall x\in \Uqg\ ,
\ee
\be\label{u-con2}
(\Delta\ot\id)\kr=\kr^{13}\kr^{23}\ ,\quad
(\id\ot\Delta)\kr=\kr^{13}\kr^{12}\ ,
\ee
\ese
where for $\Delta(x)=x_{(1)}\ot x_{(2)}$, $\Delta'(x)=x_{(2)}\ot x_{(1)}$
means the opposite comultiplication and if $\kr=\sum a_i\ot b_i$,
$\kr^{13}$ means $\sum a_i\ot1\ot b_i$, etc.

It is known \cite{KT} that
the universal $\kr$-matrix for quantum affine algebra with Hopf structure
given by \r{copr} have the following form
\be\label{PP2}
\kr_{\rm can}=\kr_{-,+}^{21}\cdot\ck\cdot\kr_{+,-}\ ,
\ee
where $A^{21}$ means the transposition of the left and right tensor space
and the element $\ck$
\be\label{WW4}
\ck =q^{-\frac{h\ot h}{2}}q^{\frac{-c\ot d -d\ot c}{2}}
\exp\left((q^{-1}-q)\sum_{n>0} \frac{n}{[2n]_q}\
a_{n}\ot a_{-n}\right)
q^{\frac{-c\ot d -d\ot c}{2}}
\ee
depends only on the Cartan and imaginary root generators.
The $\kr$-matrix \r{PP2}
belongs to the tensor product
$U_q(\widehat{\frak{b}_+})\ot
U_q(\widehat{\frak{b}_-})$, where
$U_q(\widehat{\frak{b}_+})$ is generated by $e_{m}$, $m\geq0$,
$k^{\pm1}$, $a_{n}$,  $f_{n}$, $n>0$ and
$U_q(\widehat{\frak{b}_-})$ is generated by
 $f_{n}$, $n\leq0$, $k^{\pm1}$, $a_{m}$, $e_{m}$, $m<0$.
The multiplicative expressions for the elements
$\kr_{\pm,\mp}$ are known (see \cite{KT} or formulas \r{krpm} below).

The main subject of the paper is an integral presentation for the
factors  $\kr_{\pm,\mp}$ of the universal $\kr$-matrix \r{PP2}.
Let  $d_\al$ be the following gradation operator:
\be\label{grad} [d_\al,e(z)]=e(z),\quad
[d_\al,f(z)]=-f(z),\quad [d_\al,\psi_\pm(z)]=0 \ee
and (do not confuse $\tau$ with spectral parameter)
\be\label{rtau}
\kr_{\pm,\mp}(\tau)=\tau^{-d_\al\ot1}\ \kr_{\pm,\mp}\ \tau^{ d_\al\ot1}\ ,
\qquad
\ov\kr_{\pm,\mp}(\tau)=1\ot1+\sum_{n>0}\  \kr^{(n)}_{\pm,\mp} \tau^n\ .
\ee

The following differential equations can be regarded as a main result
of the paper.
\begin{th}\label{difur}
Let $q^N\not=1$ for $N\in\ZZ\setminus\{0\}$. Then
\bse\label{I55}
\be\label{I55a}
\tau\frac{d{\kr}_{+,-}(\tau)}{d \tau}=
{\kr}_{+,-}(\tau)\cdot I_{+,-}(\tau)
 \ ,
\ee
\be\label{I55b}
\tau\frac{d{\kr}_{-,+}(\tau)}{d \tau}=
 I_{-,+}(\tau)\cdot {\kr}_{-,+}(\tau)\ ,
\ee
\ese
where
\be\label{hom-c}
I_{\pm,\mp}(\tau)=\sum_{n>0} I^{(n)}_{\pm,\mp} \tau^n
\ee
and
\be\label{I10}
I^{(n)}_{+,-}
=\frac{(-1)^n(q^{-1}-q)}{[n]_q![n-1]_q!}\oint\frac{dz}{2\pi i z}\ S_{f_0}^{n-1}
\sk{\fp(z)}\otimes
S_{e_0}^{n-1}\sk{\em(z)}\ ,\\
I^{(n)}_{-,+}
=\frac{(-1)^n(q^{-1}-q)}{[n]_q![n-1]_q!}
\oint\frac{dz}{2\pi iz}\ {S}_{f_0}^{n-1}
\sk{\fm(z)}\otimes
{S}_{e_0}^{n-1}\sk{\ep(z)}\ .
\ee
\end{th}
The screening operators $S_{e_0}$ and $S_{f_0}$ are defined
through left/right adjoint actions \r{defin} which use the standard Hopf
structure  \r{copr}:
$$S_{e_0}(x)=e_0x-k_{}xk_{}^{-1}e_0,\qquad
S_{f_0}(x)=xf_0-f_0k_{}^{-1}xk_{}.$$
It is possible also to express the action of these screening operators on the
fields $e_\pm(z)$ and $f_\pm(z)$ via the powers of the fields:
$$S^{n-1}_{e_0}\skk{e_\pm(z)}=\prod_{k=2}^n(1-q^{ 2(k-1)}) e_\pm^n(z),\qquad
S^{n-1}_{f_0}\skk{f_\pm(z)}=\prod_{k=2}^n(q^{- 2(k-1)}-1) f_\pm^n(z).$$

The differential equations \r{I55a} and \r{I55b} define the recurrence
 relations between homogeneous components $\kr_{\pm,\mp}^{(n)}$ of
 $\kr_{\pm,\mp}(\tau)$. Moreover, these equations have unique solutions
 in power series over $\tau$ with initial conditions
$\kr_{\pm,\mp}(0)=1\ot1$:
\begin{th}\label{main-theorem}
The elements
$\kr_{\pm,\mp}$ can be presented as series of multiple formal integrals
\be\label{int11}
\kr_{\pm,\mp}=1\ot1+\sum_{n>0}\ \ov\kr^{(n)}_{\pm,\mp},
\ee
\be\label{int2}
\kr^{(n)}_{\pm,\mp}=(-2\pi i)^{-n}\sum_{m=1}^n\sum_{j_1+j_2+\ldots+j_m=n}
C_{\pm}(j_1,j_2,\ldots,j_m)\times\\
\times \oint\cdots\oint\frac{dz_1}{z_1}\cdots \frac{dz_m}{z_m}\
S_{f_0}^{j_1-1}\sk{f_{\pm}(z_1)}\cdots S_{f_0}^{j_m-1}\sk{f_{\pm}(z_m)}
\ot S_{e_0}^{j_1-1}\sk{e_{\mp}(z_1)}\cdots S_{e_0}^{j_m-1}\sk{e_{\mp}(z_m)}
\ee
and
\be\label{int3}
C_+(j_1,j_2,\ldots,j_m)=\frac{(q^{-1}-q)^{m}}
{j_1(j_1+j_2)(j_1+j_2+j_3)\cdots (j_1+j_2+\cdots+j_m)}
\prod_{i=1}^m \frac{1}{[j_i]_q![j_i-1]_q!},\\
C_-(j_1,j_2,\ldots,j_m)=\frac{(q^{-1}-q)^{m}}
{j_m(j_m+j_{m-1})(j_m+j_{m-1}+j_{m-2})\cdots (j_m+j_{m-1}+\cdots+j_1)}
\prod_{i=1}^m \frac{1}{[j_i]_q![j_i-1]_q!}.
\ee
\end{th}
Applying the results of \cite{DM}, we can prove that in integral
representations of level $k>0$
  the fields $e_\pm(z)$ and $f_\pm(z)$  are annihilated by
$k+1$-th degree of screenings. So we have under the same condition on $q$

\begin{cor}\label{collorary}
Let the  $\kr$-matrix \r{PP2} acts in tensor product of integrable
 representations, one of which has level $k>0$. Then the summation
indices $j_i$,
$i=1,\ldots,m$
in \r{int2} satisfy the inequalities $1\leq j_i\leq k$.
 In particular, if one of the representations has level $1$, then the
 the  $\kr$-matrix \r{PP2} has a form
\be\label{cor1}
\kr_{\rm can}=
\exp\sk{\frac{q^{-1}-q^{}}{2\pi i}\oint \frac{dz}{z}\ e_+(z)\ot f_-(z)}\ \cdot
\ck\cdot \exp\sk{\frac{q^{-1}-q^{}}{2\pi i}\oint \frac{dz}{z}\ f_+(z)\ot
e_-(z)}\ ,
\ee
where the factor $\ck$ is given by $\r{WW4}$.
\end{cor}

\section{$U_q(\widehat{sl}_2)$ with `Drinfeld' coproduct}
\label{Drinfeld}

The goal of this section is to present different expressions for the
universal $\kr$-matrix in quantum affine algebra $\Uqdva$ with respect to
Drinfeld comultiplication. We describe it in
multiplicative form \cite{KT} and  in a form of contour integral \cite{DK}.
We also present formal integral expression for the tensor of the pairing of
two Borel subalgebras \cite{ER,EF}.

In the first subsection we describe so called Drinfeld comultiplication for
$\Uqdva$, which naturally appear in his 'new realization'. Then we study the
Hopf pairing between corresponding Borel subalgebras and remind the
multiplicative expression for the universal $\kr$-matrix, which diagonalize
this pairing. In this form the universal $\kr$-matrix can act either in tensor
product of highest weight representations or in tensor
product of lowest weight representations.

Next we remind the construction from \cite{DK} of the universal $\kr$-matrix
 in a form of contour integral. Its description requires the use of a
 natural completion of $\Uqdva$ and can act only on tensor product of
highest weight representation. Though such element cannot be factorized in
desired way, the differential equation \cite{DKP} on it will be very
important further.

Finally, we observe a presentation of \cite{EF,ER} of the tensor of the
pairing as simple formal integral. We stress certain delicate details in
this presentation: its makes sense only under specific normal ordering,
which allows to apply the formal integral to tensor product of lowest and
highest weight representations. Due to this antisymmetry, formal integral
 is difficult to use as an $\kr$-matrix with action in certain tensor
category of representations, while the factorization into product of
two cocycles makes sense.

\subsection{The universal $\kr$-matrix in a multiplicative form}

In \cite{D1} the quantum affine algebras have been constructed
by means of the quantum double construction with comultiplication
different from those given in formulas \r{copr}.
This is so called Drinfeld Hopf structure written in terms of generating
functions as follows:
\bse\label{du1}
\be\label{du1a}
\Delta^{(1)}\xp(z)=\xp(z)\otimes 1+\psi^-(zq^{\frac{c_1}{2}})
\otimes \xp(zq^{c_1})\ ,
\ee
\be\label{du1b}
\Delta^{(1)}\xm(z)=1\otimes \xm(z)+\xm(zq^{c_2})\otimes
\psi^+(zq^{\frac{c_2}{2}})\ ,
\ee
\be\label{du1c}
\Delta^{(1)}\psi^\pm(z)=\psi^\pm(zq^{\pm\frac{c_2}{2}})\otimes
\psi^\pm(zq^{\mp\frac{c_1}{2}})\ ,
\ee
\be\label{du1d}
\ant^{(1)}\sk{e(z)}=-\sk{\psi^-\sk{zq^{-\frac{c}{2}}}}^{-1}\ e(zq^{-c}),\quad
\ant^{(1)}\sk{f(z)}=-f(zq^{-c})\ \sk{\psi^+\sk{zq^{-\frac{c}{2}}}}^{-1}\ ,
\ee
\be\label{du1f}
\ant^{(1)}\sk{\psi^\pm(z)}=\sk{\psi^\pm(z)}^{-1},\quad
\ant^{(1)}(c)=-c,\quad \ant^{(1)}(d)=-d\ ,
\ee
\be\label{du1e}
\coun(c)=
\coun(d)=\coun(e(z))=\coun(f(z))=0,\quad
\coun(\psi^\pm(z))=1\ .
\ee
\ese
There exists two Drinfeld's Hopf structures.
That is why we denoted coproduct and  antipode  maps
 in \r{du1} as $\Delta^{(1)}$ and $\ant^{(1)}$. The counit maps
are the same in the all Hopf structures considered in this paper, so
we do not use for them different notations. The second
structure is given  by the formulas
\bse\label{du21}
\be\label{du21a}
\Delta^{(2)}\xp(z)=\xp(z)\otimes 1+\psi^+(zq^{-\frac{c_1}{2}})
\otimes \xp(zq^{-c_1})\ ,
\ee
\be\label{du21b}
\Delta^{(2)}\xm(z)=1\otimes \xm(z)+\xm(zq^{-c_2})\otimes
\psi^-(zq^{-\frac{c_2}{2}})\ ,
\ee
\be\label{du21d}
\ant^{(2)}\sk{e(z)}=-\sk{\psi^+\sk{zq^{\frac{c}{2}}}}^{-1}\ e(zq^{c}),\quad
\ant^{(2)}\sk{f(z)}=-f(zq^{c})\ \sk{\psi^-\sk{zq^{\frac{c}{2}}}}^{-1}
\ee
\ese
and the rest maps are the same as in \r{du1}.

To exploit the  method of quantum double construction we should
decompose the algebra under consideration into two Borel subalgebras
corresponding to the chosen Hopf structure. This decomposition for \r{du1}
is\footnote{As usual, one should
start with different Cartan elements $c$, $d$ and $\psi^\pm_{0}$
in these dual subalgebras and factorizing over simple central elements
identify these Cartan elements. In this way one can prove
 the relation $\psi^+_{0}=\sk{\psi^-_{0}}^{-1}$.}
$U_F=\{\xm(z)\,,\;\psi^+(z)\,,\;c\,,\;d\}$ and
$U_E=\{\xp(z)\,,\;\psi^-(z)\,,\;c\,,\; d\}$.
The pairing between these dual subalgebras  can be
 reconstructed from comparing the general multiplication
in quantum double  with the commutation relations \r{10} and \r{7}.
The only nontrivial pairings between generators are $\<c,d\>=\<d,c\>=1$ and
(in terms of the  formal series)
\be\label{du6}
\<\xm(w),\xp(z)
\>\,=\,\,\frac{1}{q^{-1}-q}\,\delta(z/w)
\ee
and
\be\label{du8}
\<\psi^+(w),\psi^-(z)\>= g\sk{\frac{z}{w}}=
\frac{q^2-\frac{z}{w}}{1-q^2\frac{z}{w}}=
q^2+(q^2-q^{-2})\sum_{k>0}q^{2k}\frac{z^k}{w^k}\ .
\ee

The pairing \r{du6} and \r{du8}   $U_{F}\ot U_{E}\rightarrow
\CC$ is the Hopf pairing, satisfying the property
\bn
\<k_1f,k_2e\>=\<k_1,k_2\>\<f,e\>
\label{pair}
\ed
 for any $k_1\in K_+$, $k_2\in K_-$, $f\in U_f$, $e\in U_e$.
 Here $K_\pm$ are the algebras,
generated by $q^{\pm h}$, $a_{n}$ for $n\geq 0$ and $n\leq 0$
and $U_f$ and $U_e$ are nilpotent subalgebras of
$U_{F}$ and $U_{E}$ which are generated only by modes $f_n$ and $e_n$,
$n\in\ZZ$ respectively.
The property \r{pair}
 signifies that the universal $\kr$-matrix corresponding to the
coproduct $\Delta^{(1)}$ will have factorized form
\be\label{un-fac}
\kr=\ck\cdot\ov\kr\ ,
\ee
where $\ck$ is given by \r{WW4}.

To describe the element $\ov\kr$ one should consider the restriction of the
pairing \r{du6} and \r{du8} on the dual subalgebras $U_f$  and $U_e$. It
 can be obtained inductively and have the form
\bn
\<f(z_1)\cdots f(z_n),e(w_1)\cdots e(w_n)\>=\\
=(q^{-1}-q)^{-n}\sum_{\sigma\in S_n}\prod_k
\delta\left(\frac{z_k}{w_{\sigma(k)}}\right)\prod_{{k<l \atop
\sigma(k)>\sigma(l)}} g\left(\frac{z_k}{z_l}\right),
\label{pr}
\ed
where the series $g(z)$ is given by the formula \r{du8} and
 $S_n$ is  symmetric group of permutations of $n$ elements.
The Hopf pairing between $U_e$ and $U_f$, corresponding to comultiplication
 $\Delta^{(2)}$ has analogous form (see \r{pr1} below).

The element $\ov\kr$ in the multiplicative form
reads as follows \cite{KT}
\be\label{du10}
\overline{\kr}=
\stackrel{\longrightarrow}{\prod_{n\in\ZZ}}
\exp_{q^2}\Big((q^{-1}-q)\xm_{-n}\otimes \xp_{n}\Big)\ ,
\ee
where
\be\label{q-exp}
\exp_{q^2}(x)=1+x+\frac{x^2}{1+q^2}+\frac{x^3}{(1+q^2)(1+q^2+q^4)}+\cdots+
\frac{x^n}{(n)_{q^2}!}+\cdots\ \ ,
\ee
where $(a)_{q^2}=\frac{q^{2a}-1}{q^2-1}$. In terms of the pairing the
formula \r{du10} is equivalent to the relation
$$
\<f_{-n_1}^{k_1}\cdots f_{-n_s}^{k_s},e_{m_1}^{l_1}\cdots e_{n_s}^{l_s}\>=
\prod_i\delta_{n_i,m_i}\cdot \delta_{k_i,l_i}\cdot (q^{-1}-q)^{k_i}(k_i)_{q^2}!$$
 for all $n_1<...< n_s$, $m_1<...<m_s$, which can be deduced from \r{pr} by
induction.

Note that the ordering of the $q$-exponents in \r{du10} dictates that modes
$f_n$
in the first tensor space are ordered by nonincreasing indices while $e_n$
in the second tensor space by nondecreasing order, so the element $\ov\kr$
belongs to the tensor product of
certain completion of the nilpotent subalgebras $U_f$ and $U_e$,
such that the element \r{du10} is well defined operator when it acts in the
tensor product of lowest weight representation with arbitrary one or
 in tensor product of arbitrary representation with highest weight
representation.

A distinguished feature of the pairing \r{pr} which follows from the form
\r{du10} is that it provides a way
to order any monomial of generators $e_n$
 as a sum of monomials $e_{k_1}\cdots e_{k_n}$ with nondecreasing indices
$k_i$, $k_1\leq..\leq k_n$ and
to order any monomial of generators $f_n$
 as a sum of monomials $f_{l_1}\cdots f_{l_n}$ with nonincreasing indices
$l_i$, $l_1\geq ..\geq l_n$. For the orderings in opposite directions
one should use the pairing corresponding to comultiplication $\Delta^{(2)}$:
\bn
\<e(w_1)\cdots e(w_n),f(z_1)\cdots f(z_n)\>^{(2)}
=(q^{}-q^{-1})^{-n}\sum_{\sigma\in S_n}\prod_k
\delta\left(\frac{z_k}{w_{\sigma(k)}}\right)\prod_{{k<l \atop
\sigma(k)>\sigma(l)}} g'\left(\frac{z_l}{z_k}\right),
\label{pr1}
\ed
where
\be\label{g'-series}
g'(z)=\frac{q^{-2}-z}{1-q^{-2}z}=q^{-2}+(q^{-2}-q^2)\sum_{k>0}
q^{-2k}z^k\ .
\ee
 An application of this technique is presented in the Appendix.

\subsection{Integral presentation for the element $\ov\kr'$ and
the fundamental differential equation}
It is a common place to extend  highest (lowest) weight representations of a
 Kac-Moody
algebra to a representations of bigger algebra. This algebra can be defined
as a completion with respect to a minimal topology, compatible with action
on highest weight representations \cite{DK}. Since we are interested only by
 separate
action of the algebras $U_e$ and $U_f$, these completions can be defined
quite explicitly. Namely, the completed algebras $\ov U_e^<$ and
 $\ov U_f^<$, which act in highest weight representations, are generated as
linear spaces by the series over monomials $e_{k_1}\cdots e_{k_l}$
($f_{k_1}\cdots f_{k_l}$) where $k_1\leq...\leq k_l$ and $\sum_{i}{k_i}$ is
fixed. Analogously, for the completed algebras $\ov U_e^>$ and
 $\ov U_f^>$, which act in lowest weight representations, we use the series
over monomials, ordered in opposite direction.

The following presentation of the universal $\kr$-matrix was constructed
 in \cite{DK} in a completed tensor product of $\ov U_f^<\ot \ov U_e^<$:
\bn \label{Rprime}
{\cal R}'={\cal K}\cdot\ov{\kr}',
\ed
where ${\cal K}$ coincides with \r{WW4} and
\be\label{DiKh}
\ov\kr'=
1+\sum_{n>0}\frac{1}{n!(2\pi i)^n}
\stackreb{D_n}{\oint\cdots\oint}\frac{dz_1}{z_1}\cdots \frac{dz_n}{z_{n}}\
t(z_1)\cdots t(z_n).
\ee
Here
\bn
t(z)=(q^{-1}-q)f(z)\otimes e(z)
\label{tz}
\ed
and $D_n$ is $n$-dimensional torus $|z_i|=1$ for $|q|>1$ and is the
$n$-cycle $|z_i\prod_{j=1,...,n,j\not=i}(z_i-q^2z_j)|=1$, $i=1,...,n$ for
 any $q$, such that $q^N\not=1$, $N\in\ZZ\setminus\{0\}$
and the integrand is understood as analytical continuation
 of the product $t(z_1)\cdots t(z_n)$ from the region
$|z_1|\gg|z_2|\gg...\gg|z_n|$ \cite{DKP}.

 The results of the paper are strongly based on the following
\begin{prop}\label{DKP-main}\cite{DKP}
{\rm (i)} The action of $\kr$-matrix \r{Rprime}
 in tensor product of highest weight modules is well defined and
coincides with the action of the $\kr$-matrix \r{un-fac};

{\rm (ii)} The element
 $\ov{\kr}'_{}(\tau)=\tau^{- d_\al\ot 1}\ \ov{\kr}'_{}\ \tau^{ d_\al\ot1}$
 satisfy the following differential equation:
\bn\label{Int-Eq}
\tau\frac{d\ov{\kr}'(\tau)}{d \tau}=\ov{\kr}'(\tau)\cdot
I(\tau)=
I(\tau)\cdot\ov{\kr}'(\tau).
\ed
\end{prop}
Here the generating series $I(\tau)$ also belongs to the tensor product
of the same completions $\ov U^<_f\ot\ov U^<_e$ and so is a well
defined operator acting in the tensor product of h.w.r. The
coefficients of  the formal series $I(\tau)$
\be\label{I1}
I(\tau)=\sum_{n=1}^\infty I^{(n)} \tau^n
\ee
 are the commuting quantities
\be\label{I0}
I^{(n)}=\frac{(q-q^{-1})(-1)^{n}}
{[n-1]_q![n]_q!}\oint\dzzz{}\ f^{(n)}(z)\otimes e^{(n)}(z)
\ee
constructed from the composed currents
 $f^{(n)}(z)$ and $e^{(n)}(z)$ (their definitions of the multiple residues
are given by \r{pp10})
\be\label{I2}
f^{(n)}(z)=(q^{-1}-q)^{n-1}[n]_q![n-1]_q!f(q^{2(n-1)}z)f(q^{2(n-2)}z)\cdots
 f(z)
\ ,\\
e^{(n)}(z)=(q-q^{-1})^{n-1}[n]_q![n-1]_q!e(z)e(q^{2}z)\cdots e(q^{2(n-1)}z)
\ .
\ee
The composed currents can be
obtained inductively as the residues of the products
of the preceding  currents
and by the constructions the composed currents are elements
from the completions $\ov{U}^<_{f,e}$.
The commutativity of the zero modes of two-tensors
\be\label{I3}
t^{(n)}(z)=\frac{(q-q^{-1})(-1)^{n}}
{[n-1]_q![n]_q!}\ f^{(n)}(z)\otimes e^{(n)}(z)
\ee
follows from the commutation relations
\be\label{I4}
[\T{n}(z_1),\T{m}(z_2)]=
\delta(q^{2n}z_1/z_2)\T{n+m}(z_1)-
\delta(q^{-2m}z_1/z_2)\T{n+m}(z_2)\ .
\ee
Due to this commutativity the differential equation can be easily solved
\be\label{solutionD}
\ov\kr'=\exp\sk{\sum_{n>0}\frac{I^{(n)}}{n}}=\exp\sk{(q-q^{-1})
\sum_{n>0}\frac{(-1)^{n}}
{[n-1]_q![n]_q!}\oint \dzzz{}\ f^{(n)}(z)\ot e^{(n)}(z)}.
\ee

\subsection{The pairing tensor as a  formal integral \cite{ER,EF}}

In the completions $\ov U^<_e$, $\ov U^>_e$, $\ov U^<_f$ and $\ov U^>_f$
 the defining relations
for  the total currents  can be strengthened.
\begin{lem} \label{lemtop}
{\rm (i)} In the algebras
$\ov{U}_{e}^<$ and  $\ov{U}_{f}^<$ the following identity of formal power
series take place:
\be \label{L31}
e(z_1)\cdots e(z_n)=\prod_{k<l}g\left(\frac{z_l}{z_k}\right)
e(z_n)
\cdots e(z_1),\quad
f(z_1)\cdots f(z_n)=
\prod_{k<l}g'\left(\frac{z_l}{z_k}\right)f(z_n)\ .
\cdots f(z_1)
\ee

\noindent {\rm (ii)}
 In the algebras   $\ov{U}_{f}^>$ and
$\ov{U}_{e}^>$ the following identity of formal power series take place:
\be \label{L32}
e(z_1)\cdots e(z_n)=
\prod_{k<l}g'\left(\frac{z_k}{z_l}\right)e(z_n)
\cdots e(z_1),\quad
f(z_1)\cdots f(z_n)=
\prod_{k<l}g\left(\frac{z_k}{z_l}\right)f(z_n)
\cdots f(z_1)\ ,
\ee
where $g(z)$ and $g'(z)$  are given by the series \r{du8} and \r{g'-series}
respectively.
\end{lem}

\noindent{\it Proof.}
One can prove the statement of this Lemma using pairing arguments
and explicit formula \r{pr} for the pairing.
For instance,  the pairing \r{pr} can be extended by
continuity to the pairing of $\ov U_e^<$ with $U_f$. Then
$$(q^{-1}-q)^2\left\<f(w_1)f(w_2), e(z_1)e(z_2) - g\sk{\frac{z_2}{z_1}}
e(z_2)e(z_1)\right\>=$$
$$
=\delta\sk{\frac{z_1}{w_1}}\delta\sk{\frac{z_2}{w_2}}
+g\sk{\frac{z_2}{z_1}}\delta\sk{\frac{z_2}{w_1}}\delta\sk{\frac{z_1}{w_2}}
-$$
$$
-g\sk{\frac{z_2}{z_1}}\left[
\delta\sk{\frac{z_2}{w_1}}\delta\sk{\frac{z_1}{w_2}}
+g\sk{\frac{z_1}{z_2}}
\delta\sk{\frac{z_1}{w_1}}\delta\sk{\frac{z_2}{w_2}}
\right]\equiv 0\ ,
$$
where the vanishing is due to the functional relation
$g(z)g(z^{-1})=1$. Note that  for two of four equalities in Lemma
\ref{lemtop} we should use the Hopf pairing \r{pr1} attached to the coproduct
$\Delta^{(2)}$.

Let us introduce the following formal integral
\bn\label{ttR}
\ov{\kr}''={:}\  \exp\sk{\oint\frac{dz}{2\pi  z}\ t(z)}\ {:}=
\sum_{n\geq 0}\frac{(q^{-1}-q)^n}{n!(2\pi i)^n}
{:}\oint\cdots \oint\dzz{1}\cdots \dzz{n}\
 f(z_1)\cdots f(z_n)\ot
e(z_1)\cdots e(z_n)\ {:}\ ,
\ee
 The dots ${:}\ \cdot\ {:}$
mean  the following ordering of the result of the formal integration:
we present the monomial of generators $e_n$
 as a sum of monomials $e_{k_1}\cdots e_{k_n}$ with increasing indices
$k_i$, $k_1\leq\cdots\leq k_n$ and
to order any monomial of generators $f_n$
 as a sum of monomials $f_{l_1}\cdots f_{l_n}$ with decreasing indices
$l_i$, $l_1\geq\cdots\geq l_n$.
The procedure is correctly defined only for $|q|<1$
since the resulting coefficients at ordered monomials include
 the sums of geometric progressions. The element $\ov\kr''$ belongs, by a
construction, to a completed tensor product of
$\ov{U}_f^>$ and $\ov{U}_e^<$.
 \begin{lem}\label{spariw}
 \bn
\<\ov{\kr}'', x\ot 1\>=x, \qquad \<\ov{\kr}'', 1\ot y\>=y
\label{Rxy}
\ed
  for any $x\in \ov{U}_f^>$ and $y\in \ov{U}_e^<$.
\end{lem}
Let us calculate, for instance, $\<\ov{\kr}'', e_(z_1)e(z_2)\ot 1\>$.
 Due to the formulas of the pairing, we have after integrations of
delta-functions:
$$
\<\ov\kr'', e(z_1)e(z_2)\ot 1\>= \frac{1}{2}\left(e(z_1)e(z_2)+
 g(z_2/z_1)e(z_2)e(z_1)\right)
$$
 which is equal to $e(z_1)e(z_2)$ due to Lemma \ref{lemtop}.
Due to the Lemma \ref{spariw},
we can use $\ov{\kr}''$ as a tensor of pairing
 between  $\ov{U}_f^>$ and $\ov{U}_e^<$. Moreover, we will see further, that
it
 actually coincides with  \r{du10} and thus use further for
 $\ov{\kr}''$ the same notation $\ov{\kr}$.

Nevertheless, from the definition we can define an
 action of  \r{ttR} only on tensor product of lowest weight and
highest weight
 representations; also the quantum double, including
$\ov{U}_f^>$ and $\ov{U}_e^<$ as dual Hopf subalgebras is not
 correctly defined due to divergences.
Thus we use   the expression \r{ttR} not as the universal
$\kr$-matrix, but as a tensor of pairing.

\section{Factorization of the universal $\kr$-matrix}
 \label{factor-universal}
As it was mentioned in Introduction, we want, following \cite{EF,ER}, to
decompose the factor $\ov\kr$ of the universal $\kr$-matrix \r{un-fac},
 corresponding to Drinfeld comultiplication, into a product of two
 cocycles, which can be used for restoring
of the canonical comultiplication structure
\r{copr} and corresponding $\kr$-matrix. For multiplicative presentation of
 $\ov\kr$ \r{du10} such a factorization is clear:
$\ov\kr=\kr_{+,-}\cdot\kr_{-,+}$ where
$\kr_{+,-}$ consists of the product of $q$-exponents
$\exp_{q^2}\Big((q^{-1}-q)\xm_{-n}\otimes \xp_{n}\Big)\ $
 for $n<0$ and
$\kr_{-,+}$ consists of the product of $q$-exponents
 for $n\geq 0$, such that the universal $\kr$-matrix for canonical
 comultiplication has a form
$\kr_{can}=\kr_{-,+}^{21}\cdot{\cal K}\cdot\kr_{+,-}$ \cite{TK}, that is,
\bn\label{krpm}
\kr_{\rm can}=
\stackrel{\longrightarrow}{\prod_{n\geq 0}}
\exp_{q^2}\Big((q^{-1}-q)e_{n}\otimes f_{-n}\Big)\
\cdot {\cal K}\ \cdot
\stackrel{\longleftarrow}{\prod_{n> 0}}
\exp_{q^2}\Big((q^{-1}-q)f_{n}\otimes e_{-n}\Big)\ .
\ed

Our goal is to develop general technique of factorization
applicable to a situation when multiplicative expression
for $\kr$-matrix
 is missing.
 Algebraical background for such a factorization consists of the use of
 projection operators from (current) Borel subalgebras to their
orthogonal subalgebras, developed in \cite{EF,ER}. Its survey is given in
the first subsection and is applied to the tensors from previous subsection
 in the second one.

The projections of the contour integrals $\ov\kr'$ are not well defined.
Still, we notice in the next subsection, that the projections of the
 logarithmic derivatives of $\ov\kr'$ make sense and we use this to deduce
 differential equations for the factors $\kr_{\pm,\mp}$, which
 determine them uniquely. The main technical problem reduces to the
calculation of the projections of the composed currents, entering into
differential equations for $\ov\kr'$. This is done is the last subsection with
a help of screening operators. These calculations finish the proof of
main theorems.

\subsection{The biorthogonal decompositions of Hopf algebras}
\label{background}

Let $\ca$ be a bialgebra with unit and counit,
$\ca_1$ and $\ca_2$ be two subalgebras of $\ca$
 satisfying the following conditions:

(i) Algebra $\ca$ admits a decomposition   $\ca=\ca_1\ca_2$,
that is the multiplication
map
\bn
\mu:\ca_1\ot \ca_2\to \ca
\label{mult}
\ed
is an isomorphism of linear spaces;

(ii) $\ca_1$ is left coideal, $\ca_2$ is right coideal:
\bn
\Delta(\ca_1)\subset \ca\ot \ca_1\ ,\qquad \Delta(\ca_2)\subset \ca_2\ot\ca\ .
\label{coideal}
\ed

Then the operators
$$\Pia{1}: \Pia{1}(a_1a_2)=a_1\varepsilon(a_2)\ ,
\quad
\Pia{2}: \Pia{2}(a_1a_2)=\varepsilon(a_1)a_2\ , \quad a_1\in \ca_1\ ,
\quad
a_2\in \ca_2
$$
are well defined projection operators from $\ca$ to $\ca_i$, satisfying the
 following property:
\bn
\mu(\Pia{1}\ot \Pia{2})\Delta (a)=a
\label{addition}
\ed
for any $a\in \ca$. Here  $\varepsilon$ is counit. In Sweedler
notation \r{addition} means
$$\Pia{1}(a')\Pia{2}(a'')=a\ .$$

 The correctness of definition of $\Pia{i}$ follows from the condition
(i). Denote by $\phi:\ca\to \ca$ the linear map
$\phi(a)=\Pia{1}(a')\Pia{2}(a'')$. We claim that
 $\phi$ is the map of left $\ca_1$-modules and of right $\ca_2$-modules,
 that is, $\phi(a_1 a)=a_1\phi(a)$, $\phi(a a_2)=\phi(a)a_2$.
Indeed,  for any $a_1\in \ca_1$, $a\in \ca$ we have
$$
\phi(a_1a)=\Pia{1}(a_1'a')\Pia{2}(a_1''a'') .
$$
{}From \r{coideal} we know that $a_1''\in A_1$ , so $\Pia{2}(a_1''a'')=
\varepsilon(a_1'')\Pia{2}(a'')$ and
$$
\phi(a_1a)=\Pia{1}(a_1'\varepsilon(a_1'')a')\Pia{2}(a'')=
\Pia{1}(a_1a')\Pia{2}(a'') =
a_1\Pia{1}(a')\Pia{2}(a'')=a_1\phi(a)\ .
$$
In analogous manner we prove, that $\phi(a a_2)=\phi(a)a_2$.
Noting that  $\phi(1)=1$,
we conclude that $\phi(a)=a$ for any $a\in \ca$.

Let now $\cb$ be bialgebra  dual to $\ca$ with opposite comultiplication,
 that is there exists nondegenerate Hopf pairing $\<,\>:\ca\ot \cb\to \CC$,
satisfying the conditions
$$
\<a, b_1b_2\>=\<\Delta(a),b_1\ot b_2\>\ , \qquad
 \<a_1a_2,b\>=\<a_2\ot a_1,\Delta(b)\>\ ,
$$
and $\kr=\sum a^\a\ot b_\a$ be the tensor of the pairing.
Let $\kr_i=(\Pia{i}\ot \id)\kr$.
 The addition identity \r{addition} yields the factorization
\bn
\kr=\kr_1\cdot \kr_2\ .
\label{R1R2}
\ed
Indeed, the tensor $\kr$ is uniquely characterized by one of the properties
$$
\<\kr,b\ot 1\>=b, \quad \mbox{\rm for any\,} b\in \cb,\qquad
\<\kr,1\ot a\>=a
\quad \mbox{\rm for any\,} a\in \ca\ .
$$
Let us calculate $\<\kr_1\kr_2,1\ot a\>$.
We have
$$
\<\kr_i,1\ot a \>=\<(\Pia{i}\ot \id)\kr, 1\ot a\>=\Pia{i}\<\kr, 1\ot a\>=
\Pia{i}(a)\ .
$$
Then
$$
\<\kr_1\kr_2,  1\ot a\>=
\<\kr_1,1\ot a'\>\<\kr_2,1\ot a''\>=\Pia{1}(a')\Pia{2}(a'')=a
$$
due to \r{addition}. It proves \r{R1R2}.

We can get a factorization of $\kr$, $\kr=\tkr_2\tkr_1$ if we start
 from the decomposition $\cb=\cb_1 \cb_2$ into a product of two subalgebras,
 being right and left (since the opposite comultiplication is used in the
pairing)
coideals of $\cb$: $\Delta(\cb_1)\subset \cb_1\ot \cb$,
$\Delta(\cb_2)\subset \cb\ot \cb_2$ and use $\tkr_i=(1\ot \Pib{i})\kr$, where
 $\Pib{1}(b_1b_2)=b_1\varepsilon(b_2)$, $\Pib{2}(b_1b_2)=b_1\varepsilon(b_2)$.
   The natural question arise: when these two decomposition
coincide, that is $\kr_1=\tkr_1$, $\kr_2=\tkr_2$?  We claim the following
\begin{prop}\label{biort}
The decomposition of pairing tensor induced by the decomposition of the
algebra $\ca$ coincides with the decomposition of this tensor induced by
the decomposition of the algebra $\cb$ if
$\ca_i$ and $\cb_j$ are mutually orthogonal, that is
$$
\<a_i,b_j\>=\varepsilon(a_i)\varepsilon(b_j), \qquad
\mbox{\rm for any\,} a_i\in \ca_i, b_j\in \cb_j,\quad i\not= j\ .
$$
\end{prop}

Indeed, let us compute $\<\kr_1,b\ot a\>$ and $\<\tkr_1, b\ot a\>$ for
 any $a\in \ca$, $b\in \cb$. We have
$$
\<\kr_1,b_1b_2\ot a_1a_2\>= \<a_1\varepsilon(a_2), b_1b_2\>=
\varepsilon(a_2)\<a'_1,b_1\>\<a''_1,b_2\>.
$$
We know that $a_1''\in \ca_1$, so
$$
\<\kr_1,b_1b_2\ot a_1a_2\>= \varepsilon(a_2)\varepsilon(b_2)
\<a'_1\varepsilon(a''_1),b_1\>=
\varepsilon(a_2)\varepsilon(b_2)\<a_1,b_1\>.
$$
Analogously,
$$
\<\tkr_1,b_1b_2\ot a_1a_2\>= \<a_1a_2, b_1\varepsilon(b_2)\>=
\varepsilon(b_2)\<a_1,b''_1\>\<a_2,b'_1\>=
$$
$$
=\varepsilon(a_2)\varepsilon(b_2)\<a_1, b''_1\varepsilon(b'_1)\>=
\varepsilon(a_2)\varepsilon(b_2)\<a_1,b_1\>
$$
since $b'_1\in \cb_1$. We see that $\kr_1=\tkr_1$.
The same story takes place for other pair. We call further the decompositions
$\ca=\ca_1\ca_2$,
$\cb=\cb_1\cb_2$, satisfying the condition describe above, as biorthogonal
 decompositions of $\ca$ and $\cb\equiv(\ca^*)^{op}$.

Let now $\ca$ be a Hopf algebra and element $\kr$ is considered as an
element from
square tensor of its quantum double $\cd(\ca)$.
For any biorthogonal  decomposition  the tensors $\kr_2$ and
$\left(\kr_1^{21}\right)^{-1}$, where
$\kr_i=(\Pia{i}\ot 1)\kr= (1\ot \Pib{i})\kr\in \ca_i\ot \cb_i=
(\Pia{i}\ot \Pib{i})\kr$
are two cocycles in the double $\cd(\ca)$, that is,
$$
\kr_2^{12}\cdot\left(\Delta\ot 1\right)\kr_2=
\kr_2^{23}\cdot\left(1\ot \Delta\right)\kr_2 \ ,
$$
$$
\left(\Delta'\ot 1\right)\kr_1 \cdot \kr_1^{12}=
\left(1\ot \Delta'\right)\kr_1 \cdot \kr_1^{23}\ .
$$
Indeed, both sides belong to $\ca_2\ot \cd(\ca)\ot \cb_2$. From the other hand,
we have
$$
\left(\Delta'\ot 1\right)\kr_1 \cdot \kr_1^{12}
\cdot \kr_2^{12}\cdot\left(\Delta\ot 1\right)\kr_2=
\left(1\ot \Delta'\right)\kr_1 \cdot \kr_1^{23}
\cdot \kr_2^{23}\cdot\left(1\ot \Delta\right)\kr_2
$$
 due to the properties of universal $\kr$-matrix, so the coassociator
$$
\Phi=\kr_2^{12}\cdot\left(\Delta\ot 1\right)\kr_2
\cdot\left(\kr_2^{23}\cdot\left(1\ot \Delta\right)\kr_2\right)^{-1}
$$
can be presented also as
$$
\Phi=\left(\left(\Delta'\ot 1\right)\kr_1 \cdot \kr_1^{12}\right)^{-1}
\cdot\left(1\ot \Delta'\right)\kr_1 \cdot \kr_1^{23}
$$
and thus belongs to the intersection of $\ca_2\ot \cd(\ca)\ot \cb_2$ and
$\ca_1\ot \cd(\ca)\ot \cb_1$, which means that it has a form $1\ot d\ot 1$
for
 some $d\in \cd(\ca)$. Then the pentagon identity on $\Phi$ says that
there is
 no nontrivial coassociator of such a form.

\subsection{Application to the algebra $\Uqdva$}

Let $U_{F}^+$ be a subalgebra of $U_F$, generated by all
 $a_{n}$, $f_{n}$, $n>0 $ and $q^h$;
$U_{f}^-$ be a subalgebra of $U_f\subset U_F$, generated by
 all $f_{n}$, $n\leq0$.
 We choose them as $A_1$ and $A_2$. The corresponding
 projectors will be denoted as $\Pfpd$ and $\Pfmd$.
Let also $ U_{E}^-$ be a subalgebra of $U_E$, generated by all
 $a_{n}$,  $e_{n}$, $n <0$;
 $ U_{e}^+$ be a subalgebra of $U_e\subset U_E$, generated
by all $e_{n}$, $n\geq0$.
 We choose them as $B_1$ and $B_2$. The corresponding
 projectors will be denoted as $\Pem$ and $\Pep$. That is,
\be
\label{Pfd}
\Pfpd(a_1a_2)=a_1\coun(a_2),\qquad
\Pfmd(a_1a_2)=\coun(a_1)a_2,\qquad
a_1\in U_{F}^+,\qquad a_2\in U_{f}^-,
\ee
\be
\label{Pe}
\Pem(a_1a_2)=a_1\coun(a_2),\qquad
\Pep(a_1a_2)=\coun(a_1)a_2,\qquad
a_1\in U_{E}^-,\qquad  a_2\in U_{e}^+.
\ee
It  follows from the definition of the open sets in $\ov U_F^>$ and
$\ov U_E^<$
\cite{DKP},
that
 the projections of small enough open neighborhoods  of zero are zero, which
means that the projectors can be defined also on the completed spaces:
\be\label{choi}
\Pfpd:\ov{U}^>_F\to U^+_{F},\quad
\Pfmd:\ov{U}^>_F\to U^-_{f},\quad
\Pep:\ov{U}^<_E\to U^+_{e},\quad
\Pem:\ov{U}^<_E\to U^-_{E}\ .
\ee

The subalgebras $U_{F}^+$, $U_f^-$,  $U_{E}^-$,
 $U_e^+$ and the corresponding projectors
\r{Pfd} and \r{Pe}
satisfy all the conditions of biorthogonal
 decomposition and can be applied to decomposition of the pairing tensor
$\ck\cdot\ov\kr$.
Due to the property \r{pair} of the Hopf pairing, the factorization
\r{R1R2} in this case has the form
\be\label{PP1}
\ck\cdot\ov\kr=\ck\cdot\kr_{+,-}\cdot\kr_{-,+}\ ,
\ee
where
 \be\label{PP}
\kr_{\pm,\mp}= (\Pfpmd\ot 1)\ov\kr
=(1\ot \Pemp)\ov\kr=(\Pfpmd\ot \Pemp)\ov\kr\ .
\ee
 It coincides with the natural factorization for
multiplicative expression of the universal $\kr$-matrix, mentioned
 in the beginning of this section.

Since the projections admit prolongation to completed (in opposite
 directions) subalgebras \r{choi}, the application  of
$(\Pfpmd\ot \Pemp)$ to the tensor \r{ttR} is well defined, and
 we can repeat the arguments of the previous subsection to the factorization
 of this tensor. So we have an equality
$\ov{\kr}''=\kr''_{+,-}\kr''_{-,+}$, where $\kr''_{\pm, \mp}$ belong to the
same subalgebras as $\kr_{\pm, \mp}$  and satisfy the same properties of
the pairing with the elements of dual subalgebras. So they coincide, as well
 as $\ov{\kr}''$ and $\ov{\kr}$. By this reason we skip further ${}''$
in the notations.

We need also further another pair of projections operators, connected to the
 comultiplication $\Delta^{(2)}$. Their restrictions to the algebras
 $U_f$ and $U_e$  can be  defined as follows.
Let $U_f^+$ be a subalgebra of $U_f$, generated by all $f_n$, $n> 0$
 and
 $U_e^-$ be a subalgebra of $U_e$, generated by all $e_n$, $n< 0$. We
put
\be
\label{Pf}
\Pfm(a_1a_2)=a_1\coun(a_2),\qquad
\Pfp(a_1a_2)=\coun(a_1)a_2,\qquad
a_1\in U_{f}^-,\qquad a_2\in U_{f}^+,
\ee
\be
\label{Ped}
\Pepd(a_1a_2)=a_1\coun(a_2),\qquad
\Pemd(a_1a_2)=\coun(a_1)a_2,\qquad
a_1\in U_{e}^+,\qquad a_2\in U_{e}^-.
\ee
As before, they can be prolonged to corresponding completed algebras.

The computation of the projections for the products of the currents
 can be carried out by means of the commutation relations between
 half-currents \r{h-c} $e_\pm(z)$ and $f_\pm(z)$, where we put everywhere
$z=z_1/z_2$:
 \bse\label{c-re} \bn\label{c-rea}
\epm(z_1)\epm(z_2)=g(z^{-1})\epm(z_2)\epm(z_1) +\psi(z^{-1})
\left(z^{-1}\epm^2(z_1)+\epm^2(z_2)\right),
\ed
\bn\label{c-reb}
\epm(z_1)\epm(z_2)=g'(z)\epm(z_2)\epm(z_1) +\psi'(z)
\left(\epm^2(z_1)+z\epm^2(z_2)\right),
\ed
\bn\label{c-rec}
\ep(z_1)\em(z_2)=g(z^{-1})\em(z_2)\ep(z_1) +\psi(z^{-1})
\left(z^{-1}\ep^2(z_1)+\em^2(z_2)\right),
\ed
\bn\label{c-red}
\em(z_1)\ep(z_2)=g'(z)\ep(z_2)\em(z_1) +\psi'(z)
\left(\em^2(z_1)+z\ep^2(z_2)\right),
\ed
\ese
\bse\label{c-rf}
\be\label{c-rfa}
\fpm(z_1)\fpm(z_2)=g'(z^{-1})\fpm(z_2)\fpm(z_1) +\psi'(z^{-1})
\left(\fpm^2(z_1)+z^{-1}\fpm^2(z_2)\right),
\ee
\be\label{c-rfb}
\fpm(z_1)\fpm(z_2)=g(z)\fpm(z_2)\fpm(z_1) +\psi(z)
\left(z\fpm^2(z_1)+\fpm^2(z_2)\right),
\ed
\bn\label{c-rfc}
\fp(z_1)\fm(z_2)=g'(z^{-1})\fm(z_2)\fp(z_1) +\psi'(z^{-1})
\left(\fp^2(z_1)+z^{-1}\fm^2(z_2)\right),
\ed
\bn\label{c-rfd}
\fm(z_1)\fp(z_2)=g(z)\fp(z_2)\fm(z_1) +\psi(z)
\left(z\fm^2(z_1)+\fp^2(z_2)\right).
\ee
\ese
Here
$$
g(z)=\frac{q^{2}-z}{1-q^{2}z},\quad
g'(z)=\frac{q^{-2}-z}{1-q^{-2}z},\quad
\psi(z)=\frac{1-q^{2}}{1-q^{2}z},\quad
\psi'(z)=\frac{1-q^{-2}}{1-q^{-2}z}.
$$
For the calculation of the projection of the product of the currents, say
$\Pfpd(f(z_1)\cdots f(z_n))$ we first replace the product of the currents by
the sum of the products of half-currents $f_\pm(z_i)$, using the relation
 $f(z)=f_+(z)-f_-(z)$, then move successively all the $f_-(z_i)$ to the right and
all $f_+(z_i)$ to the left. The projector $\Pfpd$ kills all the factors
 $f_-(z_i)$  which stand from the right leaving at the end the products
 of some $f_+(z_i)$.

\subsection{Differential equations for the elements $\kr_{\pm,\mp}(\tau)$}

Let us rewrite the differential equation \r{Int-Eq} for the element
$\ov\kr'\in\ov U^<_f\ot\ov U^<_e$ in the form
\be\label{Eq-rew}
I(\tau)=\sk{\ov\kr'(\tau)}^{-1}\cdot \tau
\frac{d\ov{\kr'}(\tau)}{d \tau}=
\tau\frac{d\ov{\kr'}(\tau)}{d \tau}
 \cdot \sk{\ov\kr'(\tau)}^{-1}
\ee
and act by left and right hand sides of this equality onto tensor product of
highest weight modules over $\Uqdva$. According to the Proposition
\ref{DKP-main} we can replace the element $\ov\kr'$ by its multiplicative
counterpart \r{du10} which possesses the factorization \r{PP1}.
Equations \r{Eq-rew} will have the form
\bse\label{I5}
\be\label{I5a}
\sk{{\kr}_{+,-}(\tau)}^{-1}\cdot
\tau\frac{d{\kr}_{+,-}(\tau)}{d \tau}+
\tau\frac{d{\kr}_{-,+}(\tau)}{d \tau}
\cdot\sk{{\kr}_{-,+}(\tau)}^{-1}
=
{\kr}_{-,+}(\tau)\cdot
I(\tau)
\cdot \sk{{\kr}_{-,+}(\tau)}^{-1},
\ee
\be\label{I5b}
\sk{{\kr}_{+,-}(\tau)}^{-1}\cdot
\tau\frac{d{\kr}_{+,-}(\tau)}{d \tau}+
\tau\frac{d{\kr}_{-,+}(\tau)}{d \tau}
\cdot\sk{{\kr}_{-,+}(\tau)}^{-1}
=
\sk{{\kr}_{+,-}(\tau)}^{-1}\cdot
I(\tau)
\cdot {\kr}_{+,-}(\tau)\ .
\ee
\ese

Let us apply the projections
$\Pfp\ot\Pem$ to the equality \r{I5a} and
$\Pfm\ot\Pep$ to the equation \r{I5b}.
Let us consider \r{I5a}. It is clear that the projection $\Pfp\ot\Pem$
kills the second term in the l.h.s. of this equality and the first term
is stable under this projection.
Let us apply the same projections to the r.h.s. of \r{I5a}.
As we already mentioned the generating series $I(\tau)$ contain the elements
which belong to the tensor product
$\ov U^<_f\ot\ov U^<_e$ (see details in \cite{DKP}).
According to this completions application of the projection
$\Pfp$ and $\Pem$ means the following.
One should move all generators which belong to subalgebra $\Ufm$ to the left
in the first tensor space of the r.h.s. of \r{I5a} and kill all the
terms which contain these sort of generators on the left. Analogously, for
the projection $\Pem$ in
the second tensor space of the r.h.s. of \r{I5a} we will move all generators
which belong to $\Uep$ to the right and kill all the terms where such
generators survive on the right hand side.
We conclude that
$$ (\Pfpm\ot\Pemp){\kr}_{\mp,\pm}(\tau)=1\ot 1$$
and
\be\label{I8}
(\Pfpm\ot\Pemp)\sk{{\kr}_{\mp,\pm}(\tau)\cdot
I(\tau)
\cdot \sk{{\kr}_{\mp,\pm}(\tau)}^{-1}}=
(\Pfpm\ot\Pemp)
I(\tau)\equiv I_{\pm,\mp}(\tau)\ .
\ee
 The substitution of \r{I8} to \r{I5a} and \r{I5b} proves the following
\begin{prop}
\label{system}
The projections
${\kr}_{\pm,\mp}(\tau)=(\Pfpmd\ot \Pemp)\ov\kr(\tau) $ of the universal $\kr$
 matrix $\kr(\tau)$ (see \r{rtau}) satisfy the following differential
equations:
\be\label{I55c}
\tau\sk{{\kr}_{+,-}(\tau)}^{-1}\cdot\frac{d{\kr}_{+,-}(\tau)}{d
 \tau}=
(\Pfm\ot\Pep)
I(\tau)
 \ ,
\ee
\be\label{I55d}
\tau\frac{d{\kr}_{-,+}(\tau)}{d \tau}\sk{{\kr}_{-,+}
(\tau)}^{-1}\cdot=
(\Pfp\ot\Pem)
I(\tau)\ ,
\ee
where $I(\tau)$ is given by \r{I1} and \r{I0}.
  \end{prop}

One can see that the differential equations  \r{I55c} and \r{I55d} are
  equivalent to the following recurrence relations:
\bse\label{rec11} \be\label{r-cal} n\kr^{(n)}_{+,-}
=\kr^{(n-1)}_{+,-}I_{+,-}^{(1)}+
\kr^{(n-2)}_{+,-}I_{+,-}^{(2)}+\cdots+
\kr^{(1)}_{+,-}I_{+,-}^{(n-1)}+I_{+,-}^{(n)}\ ,
\ee
\be\label{l-cal}
n\kr^{(n)}_{-,+}
=
I_{-,+}^{(1)}\kr^{(n-1)}_{-,+}+
I_{-,+}^{(2)}\kr^{(n-2)}_{-,+}
+\cdots+
I_{-,+}^{(n-1)}\kr^{(1)}_{-,+}
+ I_{-,+}^{(n)}\ .
\ee
\ese
The system of the recurrent relations \r{rec11} with noncommutative
coefficients $I^{(n)}_{\pm,\mp}$ have formal solution
\be\label{for-sol}
\kr^{(n)}_{\pm,\mp}=\sum_{m=1}^n\sum_{j_1+j_2+\ldots+j_m=n}
{\cal C}_{\pm}(j_1,j_2,\ldots,j_m)
I^{(j_1)}_{\pm,\mp} I^{(j_2)}_{\pm,\mp}\cdots I^{(j_m)}_{\pm,\mp}\ ,
\ee
where
\be\label{ord-exp}
{\cal C}_+(j_1,j_2,\ldots,j_m)=\frac{1}
{j_1(j_1+j_2)(j_1+j_2+j_3)\cdots (j_1+j_2+\cdots+j_m)}
\ ,\\
{\cal C}_-(j_1,j_2,\ldots,j_m)=\frac{1}
{j_m(j_m+j_{m-1})(j_m+j_{m-1}+j_{m-2})\cdots (j_m+j_{m-1}+\cdots+j_1)}
\ .
\ee
So we reduced the proof of the Theorem 1 and of  Theorem 2
to the calculation of the projections of the  currents \r{I8}
onto subalgebras $U_{e,f}^\pm$.
We  solve this  problem  in the next subsection.

\subsection{Projections of composed currents and
screening operators}

To calculate the projections to the subalgebras $U^\pm_{e,f}$
 from the composed currents we need the following recurrent definitions of
these currents
\bse\label{pp10}
\be\label{pp10a}
e^{(n)}(z)=
-\res{w=z}e^{(n-1)}(zq^{2})e(w)\frac{dw}{w}=
-\oint_{w\ {\rm around}\ z} \dw\
e^{(n-1)}(zq^{2})e(w)\ ,
\ee
\be\label{pp10b}
e^{(n)}(z)=
\res{w=zq^{2(n-1)}}e(w)e^{(n-1)}(z)\frac{dw}{w}
=
\oint_{w\ {\rm around}\ zq^{2(n-1)}} \dw\
e(w)e^{(n-1)}(z)
\ ,
\ee
\be\label{pp10d}
f^{(n)}(z)=
-\res{w=zq^{2(n-1)}}f^{(n-1)}(z)f(w)\frac{dw}{w}
=
-\oint_{w\ {\rm around}\ zq^{2(n-1)}} \dw\
f^{(n-1)}(z)f(w)
\ ,
\ee
\be\label{pp10c}
f^{(n)}(z)=
\res{w=z}f(w)f^{(n-1)}(zq^{2})\frac{dw}{w}=
\oint_{w\ {\rm around}\ z} \dw\
f(w)f^{(n-1)}(zq^{2})\ .
\ee
\ese
Let us define the  screening operators $S_{e_0}$, $\tilde{S}_{e_0}$,
 $S_{f_0}$, $\tilde{S}_{f_0}$, which act as the following
 $q$-commutators in the algebra
$\Uqdva$:
$$S_{e_0}(x)=e_0x-k_{}xk_{}^{-1}e_0,\qquad
S_{f_0}(x)=xf_0-f_0k_{}^{}xk_{}^{-1},$$
$$\tilde{S}_{e_0}(x)=xe_0-e_0k_{}^{-1}xk_{},\qquad
\tilde{S}_{f_0}(x)=f_0x-k_{}^{-1}xk_{}f_0.$$
  The screening operators coincide with adjoint action of the elements
  $e_0$ and $f_0$ with respect to comultiplication \r{copr}:
\be\label{defin}
S_{e_0}(x)=e_0'\cdot x \cdot \ant\sk{e_0''},
\qquad  S_{f_0}(x)=\ant\sk{f_0'}\cdot x\cdot f_0'',\\
\tilde{S}_{e_0}(x)=\ant^{-1}\sk{e''_0}\cdot x \cdot e_0',
\qquad  \tilde{S}_{e_0}(x)=
f_0''\cdot x\cdot \ant^{-1}\sk{f_0'},
\ee
and are connected via the conjugation by $k$:
\be
\label{conjug}
\tilde{S}_{e_0}(x)=-q^{2}k^{-1}S_{e_0}(x)k,\qquad
\tilde{S}_{f_0}(x)=-q^{-2}k^{-1}S_{f_0}(x)k.
\ee
The application of the screening operators to a factorization problem
 is based in their compatibility with projection operators.
\begin{lem}\label{commutativity}
{\rm (i)} Subalgebras $U_e^\pm$  are invariant with respect
 to the screening operators $S_{e_0}$ and $\tilde{S}_{e_0}$;
 subalgebras $U_f^\pm$  are invariant with respect
 to the screening operators $S_{f_0}$ and $\tilde{S}_{f_0}$;

{\rm (ii)} The screening operators
$S_{e_0}$ and $\tilde{S}_{e_0}$
commute with the projectors $P_e^\pm$;
the screening operators
$S_{f_0}$ and $\tilde{S}_{f_0}$
commute with the projectors $P_f^\pm$:
\be
\label{per1}
P_e^\pm S_{e_0}(x)=S_{e_0}P_e^\pm(x),\qquad
P_e^\pm \tilde{S}_{e_0}(x)=\tilde{S}_{e_0}P_e^\pm(x),
\qquad x\in U_e\ ,
\ee
\be
\label{per2}
P_f^\pm S_{f_0}(x)=S_{f_0}P_f^\pm(x),\qquad
P_f^\pm \tilde{S}_{f_0}(x)=\tilde{S}_{f_0}P_f^\pm(x),
\qquad x\in U_f\ .
\ee
\end{lem}
The proof of statement (i) of the Lemma consists of a short calculation
based on the use of \r{1} and \r{2}. The statement (ii) follows from (i)
 together with a remark, that $\varepsilon S(x)=0$ for any screening
$S$ in consideration.

The main result of this subsection is the following
\begin{prop}\label{composed-projections}
The projections of the  currents
$e^{(n)}(z)$ and  $f^{(n)}(z)$
onto subalgebras
$U^\pm_{e,f}$ are given by the formulas
\bse\label{pp15}
\be\label{pp15a}
\Pep\sk{e^{(n)}(z)}= \tilde{S}^{n-1}_{e_0}\sk{e_+(zq^{2(n-1)})},\quad
\Pem\sk{e^{(n)}(z)}= - S^{n-1}_{e_0}\sk{e_-(z)}\ ,
\ee
\be\label{pp15b}
\Pfp\sk{f^{(n)}(z)}= S^{n-1}_{f_0}\sk{f_+(z)},\quad
\Pfm\sk{f^{(n)}(z)}= - \tilde{S}^{n-1}_{f_0}\sk{f_-(zq^{2(n-1)})}\ .
\ee
\ese
\end{prop}

\noindent
{\it Proof.}
One can prove by induction (see \cite{DKP} for details) that the currents
$e^{(n)}(z)$ and  $f^{(n)}(z)$ satisfy the following quadratic relations:
\bse\label{pp11}
\be\label{pp11a}
(w-zq^{2(n-2)})(w-zq^{2(n-1)}) e(w)e^{(n-1)}(z)=q^{2(n-1)}
(w-zq^{-2})(w-z)e^{(n-1)}(z)e(w)\ ,
\ee
\be\label{pp11b}
(w-zq^{2(n-2)})(w-zq^{2(n-1)}) f^{(n-1)}(z)f(w)=q^{2(n-1)}
(w-zq^{-2})(w-z)f(w)f^{(n-1)}(z)\ ,
\ee
\ese
Moreover,
the product $e(w)e^{(n-1)}(z)$ has unique simple pole at
$w=q^{2(n-1)}z$;
the product $e(w)e^{(n-1)}(z)$ has unique simple zero at
$w=q^{-2}z$;
the product $f(w)f^{(n-1)}(z)$ has unique simple pole at
$w=q^{-2}z$;
and the product $f(w)f^{(n-1)}(z)$ has unique simple zero at
$w=q^{2(n-1)}z$.

It means that the residues in \r{pp10a}-\r{pp10c} can be presented as
 the following formal integrals:
\bse\label{pp12}
\be\label{pp12a}
e^{(n)}(z)=
\oint\frac{dw}{2\pi iw}\left( e^{(n-1)}(zq^2)e(w)-
q^{-2(n-1)}
 e(w) e^{(n-1)}(zq^2)\al_n\sk{\frac{z}{w};q}\right)=
\ee
\be\label{pp12b}
=\oint \frac{dw}{2\pi iw}\left(e(w)e^{(n-1)}(z)
-q^{-2(n-1)}
 e^{(n-1)}(z)e(w)
\b_n\sk{\frac{w}{z};q}\right)\ ,
\ee
\be\label{pp12c}
f^{(n)}(z)=
\oint\frac{dw}{2\pi iw}\left(f^{(n-1)}(z)f(w)
-q^{2(n-1)} f(w)f^{(n-1)}(z)
\b_n\sk{\frac{z}{w};q^{-1}}\right)=
\ee
\be\label{pp12d}
=
\oint\frac{dw}{2\pi iw}\left( f(w)f^{(n-1)}(zq^2)-q^{2(n-1)}
 f^{(n-1)}(zq^2)f(w)
\al_n\sk{\frac{w}{z};q^{-1}}\right) ,
\ee
\ese
where
\be\label{ccff}
\al_n(x;q)=\frac{\sk{1-q^{2(n-1) }x} \sk{1-q^{2n}x}}
{\sk{1-x} \sk{1-q^{2}x}},\quad \b_n(x;q)=\frac{\sk{1-q^{2}x} \sk{1-x}}
{\sk{1-q^{-2(n-2)}x } \sk{1-q^{-2(n-1)}x}}\ .
\ee
The r.h.s. of \r{pp12a}--\r{pp12d} can be presented as total
integrals of left/right
 adjoint actions of the currents  $e(w)$ and $f(w)$ with respect
to coproduct $\Delta^{(1)}$. For instance, the relation \r{pp12a} we
 can read as
$$e^{(n)}(z)=\oint\frac{dw}{2\pi iw}ad_{e(w)}X=
\oint \frac{dw}{2\pi iw}\ant^{-1}\sk{e''(w)}\ X\ e'(w)\ , $$
where $X=e^{(n-1)}(zq^2)$.

We can rewirite \r{pp12a}--\r{pp12d} as
\bse\label{pp13}
\be\label{pp13a}
e^{(n)}(z)=
e^{(n-1)}(zq^2)e_0-q^{-2(n-1)}e_0 e^{(n-1)}(zq^2)
+\sum_{k<0}\a_{n,k}(q)e_k e^{(n-1)}(zq^2)  z^{-k}\ ,
\ee
\be\label{pp13b}
e^{(n)}(z)=
e_{0}e^{(n-1)}(z)- q^{-2(n-1)} e^{(n-1)}(z)  e_{0}
+\sum_{k>0} \b_{n,k}(q) e^{(n-1)}(z) e_k z^{-k}\ ,
\ee
\be\label{pp13c}
f^{(n)}(z)=
 f^{(n-1)}(z) f_0 - q^{2(n-1)} f_0 f^{(n-1)}(z)
+\sum_{k<0 }\b_{n,k}(q^{-1}) f_k f^{(n-1)}(z) z^{-k}\ ,
\ee
\be\label{pp13d}
f^{(n)}(z)=
f_0f^{(n-1)}(zq^2)-q^{2(n-1)}  f^{(n-1)}(zq^2)f_0
+\sum_{k>0}\a_{n,k}(q^{-1}) f^{(n-1)}(zq^2) f_k z^{-k}\ ,
\ee
\ese
where $\a_{n,k}(q)$ and $\b_{n,k}(q)$ are coefficients of the expansion
of the rational functions $\al_n(x;q)$ and $\b_n(x;q)$
into series with respect to $x$.

Due to the definitions\r{Pe} and  \r{Pf}  of the projection operators,
 the sums in the right hand sides
of \r{pp13} disappear under corresponding
projections
and we obtain
\bse\label{pp14}
\be\label{pp14a} \Pep\sk{e^{(n)}(z)}= 
\Pep\sk{e^{(n-1)}(zq^2)e_0-
{q^{-2(n-1)}} e_0 e^{(n-1)}(zq^2)}=
\Pep\tilde{S}_{e_0}e^{(n-1)}(zq^2)
\ ,
\ee
\be\label{pp14b}
\Pem\sk{e^{(n)}(z)}=
\Pem\sk{ e_0e^{(n-1)}(z)-{q^{2(n-1)}}e^{(n-1)}(z)e_0}=
\Pem{S}_{e_0}e^{(n-1)}(z)
\ ,
\ee
\be\label{pp14c}
\Pfp\sk{f^{(n)}(z)}=
\Pfp\sk{f^{(n-1)}(z)f_0-{q^{-2(n-1)}}f_0f^{(n-1)}(z)}=
\Pfp{S}_{f_0}f^{(n-1)}(z)
\ ,
\ee
\be\label{pp14d}
\Pfm\sk{f^{(n)}(z)}=
 \Pfm\sk{f_0f^{(n-1)}(zq^2)-{q^{2(n-1)}}f^{(n-1)}(zq^2) f_0}=
\Pfm\tilde{S}_{f_0}f^{(n-1)}(zq^2)
\ .
\ee
\ese
Iteration of the formulas \r{pp14} together
 with \r{per1}, \r{per2}, proves the
Proposition \ref{composed-projections}.
 The additional minus  in \r{pp15}
appear due to the definitions \r{h-c} $\Pem(e(z))=-e_-(z)$ and
$\Pfm(f(z))=-f_-(z)$. Note that the coefficients in front of
$e^{(n-1)}(z)e_0$ and $f_0f^{(n-1)}(z)$ in \r{pp14b} and \r{pp14c}
are changed in contrast to \r{pp13b} and \r{pp13c} in order
to have possibility to apply the statement of the Lemma
\ref{commutativity}.

 The combination of Propositions \ref{system} and
 \ref{composed-projections} complete the proof of Theorems 1 and 2.
  Note that in its formulations the screenings $\tilde{S}_{e_0}$ and
 $\tilde{S}_{f_0}$ do not appear. We can avoid their use since
 $$\tilde{S}_{f_0}^k(f_+(z))\ot \tilde{S}_{e_0}^k(e_-(z))=
 {S}_{f_0}^k(f_+(z))\ot {S}_{e_0}^k(e_-(z))$$
 because of \r{conjug}. The Corollary 1 to Theorem 2 is also a direct
consequence of the Proposition \ref{composed-projections}
due to \cite{DM}.

The adjoint action of the screening operators
onto half-currents can be expressed through
the powers of these half-currents. We have the following
\begin{lem}\label{lem-powers}
\be\label{powers}
\tilde S^{n-1}_{e_0}\skk{e_+(z)}=\prod_{k=2}^n(1-q^{-2(k-1)}) e_+^n(z),\quad
S^{n-1}_{e_0}\skk{e_-(z)}=\prod_{k=2}^n(1-q^{2(k-1)}) e_-^n(z)\ ,\\
S^{n-1}_{f_0}\skk{f_+(z)}=\prod_{k=2}^n(q^{-2(k-1)}-1) f_+^n(z),\quad
\tilde S^{n-1}_{f_0}\skk{f_-(z)}=\prod_{k=2}^n(q^{2(k-1)}-1) f_-^n(z)\ .
\ee
\end{lem}

\noindent {\it Proof.} Let us first equality in \r{powers} since the rest
are analogous. The proof is by induction over $n$.
For $n=2$ the identity
$$
\tilde S_{e_0}\sk{e_+(z)}=[e_0,e_+(z)]_{q^{-2}}=e_0 e_+(z)-q^{-2} e_+(z) e_0=
(1-q^{-2})e^2_+(z)\ .
$$
follows from the commutation relations \r{1}. Suppose that the identity
$\tilde S^{m-1}_{e_0}\skk{e_+(z)}=\prod_{k=2}^m(1-q^{-2(k-1)}) e_+^m(z)$
is valid for $m=2,\ldots,n-1$. Then we calculate
$$
\tilde S_{e_0}^n\sk{e_+(z)}=\prod_{k=2}^{n-1}(1-q^{-2(k-1)})\tilde S_{e_0}
\sk{e_+^{n-1}(z)}=\prod_{k=2}^{n}(1-q^{-2(k-1)})
e_+^{n}(z)
$$
which prove the Lemma.

This Lemma allows one to write down the integral formulas for the
elements $\kr^{(n)}_{\pm,\mp}$ in the form of the formal integrals from
two-tensors constructed from powers of half-currents:
\be\label{int222}
\kr^{(n)}_{\pm,\mp}=(-2\pi i)^{-n}\sum_{m=1}^n\sum_{j_1+j_2+\ldots+j_m=n}
\tilde C_{\pm}(j_1,j_2,\ldots,j_m)\times\\
\times \oint\cdots\oint\frac{dz_1}{z_1}\cdots \frac{dz_m}{z_m}\
f^{j_1}_{\pm}(z_1)\cdots f^{j_m}_{\pm}(z_m)
\ot e^{j_1}_{\mp}(z_1)\cdots e^{j_m}_{\mp}(z_m)
\ee
and
\be\label{int333}
\tilde C_+(j_1,j_2,\ldots,j_m)=\frac{(q^{-1}-q)^{2n-m}}
{j_1(j_1+j_2)(j_1+j_2+j_3)\cdots (j_1+j_2+\cdots+j_m)}
\prod_{i=1}^m \frac{1}{[j_i]},\\
\tilde C_-(j_1,j_2,\ldots,j_m)=\frac{(q^{-1}-q)^{2n-m}}
{j_m(j_m+j_{m-1})(j_m+j_{m-1}+j_{m-2})\cdots (j_m+j_{m-1}+\cdots+j_1)}
\prod_{i=1}^m \frac{1}{[j_i]}.
\ee
Note, that the presentation \r{int222} is specific for the case under
consideration and in general situation of $\Uqg$
the only possibility is to use the screening operators.

\section{Factorization of the formal pairing tensor}
\label{factor-pairing}
 The recurrence relations and corresponding differential
 equations for the factorized part of $\kr$-matrix
 can be also deduced from the factorization of the pairing
 tensor in a form of formal integral. We will see that in this approach
  they appear in a different form. Its equivalence to the results
 above reduces to certain combinatorial identity.
We will give the proof of this identity, proposed by A.~Okounkov.

\subsection{Another form of the differential equations and combinatorial
identity}

We would like to factorize the element \r{ttR}
\be\label{ttR1}
 {\ov\kr}={:}\ \exp\oint t(z)\dz\ {:}=
  \sum_{n\geq 0}\frac{1}{n!}
{:}\oint\cdots \oint\dzz{1}\cdots \dzz{n}\
 t(z_1)\cdots t(z_n)\ {:}
\ee
into a product
\be\label{fact}
\ov\kr={\kr}_{+,-}\cdot {\kr}_{-,+}\ .
\ee
Here, as before,
$t(z)=(q^{-1}-q)f(z)\ot e(z)\ .$
 If we use the notation $(a\ot b)_{\pm,\mp}$ for
$\Pfpmd(a)\ot \Pemp(b)$, then
 the relation \r{fact} means that
$$
\kr_{\pm,\mp}= \sum_{n\geq 0}\frac{1}{n!}
\oint\cdots \oint \dzz{1}\cdots \dzz{n}\
 \sk{t(z_1)\cdots t(z_n)}_{\pm,\mp}\ .
$$
Let
$$
\ov\kr^{(n)}= \frac{1}{n!}
:\oint\cdots \oint \dzz{1}\cdots \dzz{n}\
 t(z_1)\cdots t(z_n):
$$
and
$$
\kr^{(n)}_{\pm,\mp}= \frac{1}{n!}
\oint\cdots \oint \dzz{1}\cdots \dzz{n}\
 \sk{t(z_1)\cdots t(z_n)}_{\pm,\mp}\ .
$$
Then the factorization \r{fact} is equivalent to
\beq\label{ind}
\ov\kr^{(n)}=\sum_{0\leq l\leq n}
\kr^{(l)}_{+,-}\kr^{(n-l)}_{-,+}\; .
\eeq
We can calculate, for example,
$\kr^{(n)}_{+,-}$ using \r{ind} by induction over $n$:
\be\label{rec}
\kr^{(n)}_{+,-}=\frac{1}{n}
\left(:\ov\kr^{(n-1)}\oint t(z)\dzz{} :\right)_{+,-}
=\frac{1}{n}\sum_{0\leq l\leq n-1}
\kr^{(n-1-l)}_{+,-}
\sk{:\kr^{(l)}_{-,+}\oint t(z)\dzz{}:}_{+,-}
\; ,
\ee
 Denote by
 \be
 \tilde{I}_{+,-}^{(l)}=
\sk{:\kr^{(l)}_{-,+}\oint t(z)\dzz{}:}_{+,-},\qquad
 \tilde{I}_{-,+}^{(l)}=
\sk{:\oint t(z)\dzz{}\kr^{(l)}_{+,-}:}_{-,+},
 \label{tildaI}
 \ee
 and
$$
\tilde{I}_{\pm,\mp}(\tau)=\sum_{n>0} \tilde{I}^{(n)}_{\pm,\mp} \tau^n
$$
With these notations we
 can organize the recurrence  relations \r{rec}
 as  the following differential
 equations.
 \begin{prop}\label{drugoidifur}
 The following differential equations follow from the factorization
 \r{fact}:
\be\label{I55e}
\tau\frac{\partial{\kr}_{+,-}(\tau)}{\partial \tau}=
{\kr}_{+,-}(\tau)\cdot \tilde{I}_{+,-}(\tau)
 \ ,
\ee
\be\label{I55f}
\tau\frac{\partial{\kr}_{-,+}(\tau)}{\partial \tau}=
 \tilde{I}_{-,+}(\tau)\cdot {\kr}_{-,+}(\tau)\ ,
\ee
 \end{prop}
Due to the uniqueness of logarithmic derivative,
 we should have the identifications:
 \bn
 \label{twoI}
 \tilde{I}_{\pm,\mp}(\tau)=I_{\pm,\mp}(\tau).
 \ed
 The calculation of the integrals \r{tildaI} is very nontrivial
  technical problem. Nevertheless, we can get some profit comparing
 \r{I55e}, \r{I55f} with
 \r{I55a}, \r{I55b}. Comparing the expressions \r{tildaI} with \r{I10}
and \r{powers} we see that the equality \r{twoI} can be considered as an
effective way to calculate the integrals \r{tildaI}.

Using pairing arguments we will demonstrate in the Appendix
that more general then \r{twoI} identities
\bse\label{con3}
\be\label{con3a}
\oint\cdots\oint \dzz{1}\cdots\dzz{k}\
\sk{\skk{f(z_1)\cdots f(z_k)}_- \fp(z_{k+1}) }_+\ot
\sk{\skk{e(z_1)\cdots e(z_k)}_+ \em(z_{k+1}) }_-=\\
=
k!\ \frac{(1-q^2)^k}{(k+1)_{q^2}}\ \fp^{k+1}(z_{k+1})\ot\em^{k+1}(z_{k+1})
\ ,
\ee
\be\label{con3b}
\oint\cdots\oint \dzz{1}\cdots\dzz{k}\
\sk{\fm(z_{k+1}) \skk{f(z_k)\cdots f(z_1)}_+
}_-\ot
\sk{\ep(z_{k+1}) \skk{e(z_k)\cdots e(z_1)}_+
}_+=\\
=
k!\ \frac{(1-q^{2})^k}{(k+1)_{q^2}}\
\fm^{k+1}(z_{k+1})\ot\ep^{k+1}(z_{k+1})\ .
\ee
\ese
(as well as \r{twoI})
are equivalent to certain combinatorial identities. In particular,
the equality  \r{con3a} is equivalent
for $|q|<1$ to
\be\label{l-stat}
\frac{1}{(n)_{q^2}!(n+1)_{q^2}!}= (1-q^{2})^n
\sum_{0\leq \la_n\leq\ldots\leq \la_1}
C_{\{\la\}}(q)q^{4(\la_1+2\la_{2}+\cdots+n\la_n)}\ ,
\ee
where the constants $C_{\{\la\}}(q)$
are parameterized by the partition $\{\la_m\}$ of natural number $n$
\be\label{combinat1}
C_{\{\la\}}(q)=\prod \frac{1}{(\la'_i-\la'_{i+1})_{q^{2}}!}\ ,
\ee
and $\{\la'_j\}$ is dual to $\{\la_k\}$ partition: $\la'_j=\# k$, such
that $\la_k\geq j$. In \r{combinat1} we assume that $(0)_q!\equiv1$.
Alternatively, these constants $C_{\{\la\}}(q)$
can be defined as follows. Let $m:\{\la_1,\ldots,\la_n\}\to
\{m_1,\ldots,m_k\}$, $1\leq k\leq n$
 be a map described by the following rule
$$
\la_{1},\ldots, \la_{n}\to
\la_{1},\ldots, \la_{m_1},\la_{m_1+1},\ldots, \la_{m_1+m_2},\ldots,
\la_{m_1+\cdots+m_{k-1}+1},\ldots, \la_{m_1+\cdots+m_k}
$$
such that
$
\la_{1}=\ldots=\la_{m_1}>\la_{m_1+1}=\ldots=\la_{m_1+m_2}>\ldots>
\la_{m_1+\cdots+m_{k-1}+1}=\ldots=\la_{m_1+\cdots+m_k}\ .
$
Then the coefficients $C_{\{\la\}}(q)$ are given by
$
\prod_{i=1}^k \frac{1}{(m_i)_{q^2}!}\ .
$

A.~Okounkov proposed an independent proof of this identity, based on the
specialization formula for Macdonald polynomials \cite{Mac},
proved in full generality by I.~Cherednik \cite{Ch}.

Consider the Cauchy identity for Macdonald polynomials:
\be\label{Koshi}
\prod_{i,j} \frac{(tx_i y_j)_\infty}{(x_i y_j)_\infty} =
\sum_\la P_\la(x;q,t)\, Q_\la(x;q,t)\,,
\ee
where $(a)_\infty = (1-a)(1-qa)(1-q^2a)\cdots$,
where $(a)_\infty = (1-a)(1-qa)(1-q^2a)\cdots$, $P_\la$ are
Macdonald polynomials, and  $Q_\la$ are dual Macdonald polynomials,
 corresponding to  Young diagram $\la$.
Set
$$
x=(t,t^2,t^3,\dots,t^n,0,0,\dots)\,,
\quad y=(t,t^2,\dots,t^n,t^{n+1},t^{n+2},\dots)
$$
and
$$
q=0\,,
$$
which specializes Macdonald polynomials to Hall-Littlewood polynomials.
Then the left-hand side of  \r{Koshi} becomes
\be\label{levaya}
\frac{1}{(1-t)^n(n+1)_t!}\ ,
\ee
and due to the formula
  (VI.$6.11'$) from the book \cite{Mac}
for evaluating a Macdonald polynomial at a point
 $(t,t^2,\dots,t^n)$
\be\label{Mac1}
P_\la(t,\dots,t^n;0,t)=t^{n(\la)+|\la|} \prod_{i=1}^{\ell(\la)}
\frac{1-t^{n-i+1}}{1-t^{\la'_{\la_i}-i+1}} \,,
\ee
and also
\be\label{Mac2}
Q_\la(t,\dots,t^{n+1},t^{n+2},\dots;0,t)=t^{n(\la)+|\la|} \,,
\ee
the right hand side of \r{Koshi} becomes
\be\label{pravaya}
(n)_t!
\sum_{\ell(\la)\le n} \frac{t^{2\sum i\la_i}}
{\prod_{k\ge 0} [\la'_{k}-\la'_{k+1}]}\ .
\ee
In the formulas \r{Mac1} and \r{Mac2}  $n(\la)=\sum_i (i-1)\la_i$,
$|\la| = \sum_i \,\la_i$
and $\ell(\la)=\la'_1$ is the number of rows of the diagram  $\la$.
Comparing \r{levaya}  and \r{pravaya} at $t=q^2$ we obtain
\r{l-stat}.

\subsection{Some examples of calculations
and vanishing of the cross terms in the integrals}

Let us demonstrate how
 the factorization works in the simplest term $\ov\kr^{(2)}$.
It has the form
\be\label{ex1}
\frac{(q^{-1}-q)^2}{(2\pi i)^2}\frac{1}{2}\oint\oint\underline{dz_1dz_2}
(f_+(z_1)-f_-(z_1))(f_+(z_2)-f_-(z_2))
\ot(e_+(z_1)-e_-(z_1))(e_+(z_2)-e_-(z_2))
\ee
The action of projection operators
$\Pfpd\ot\Pem$  described above results in the following
formula:
\be\label{ex2}
\frac{(q^{-1}-q)^2}{(2\pi i)^2}
\frac{1}{2}\oint\oint \dzz{1}\dzz{2}
\left(
\fp(z_1)\fp(z_2)-\psi\bigl(\frac{z_1}{z_2}\bigr)\fp^2(z_2)\right)\ot
\left(
\em(z_1)\em(z_2)-\psi\bigl(\frac{z_2}{z_1}\bigr)\em^2(z_2)\right).
\ee
We
can open the brackets and look to four summands.  First, we leave
the regular term as it is.  The integrals
\be\label{cros}
\frac{1}{2\pi i}\oint
\dzz{1}\em(z_1)\psi\sk{\frac{z_1}{z_2}}\qquad
\hbox{and}\qquad \frac{1}{2\pi i}\oint
\dzz{1}\fp(z_1)\psi\sk{(\frac{z_2}{z_1}}
\ee
vanish due to
 definitions (or analytical properties of $\fp(z)$ and $\em(z)$).
  For the last
term we need to calculate the integral
$$
\frac{1}{2\pi i}\oint \dzz{1}
\psi\sk{\frac{z_1}{z_2}}\psi\sk{\frac{z_2}{z_1}}
\frac{1-q^2}{1+q^2}
$$
to obtain
$$
\kr^{(2)}_{+,-}=
\frac{(q^{-1}-q)^2}{2}\sk{\oint\oint
\dzzz{1}\dzzz{2} \fp(z_1)\fp(z_2)\ot \em(z_1)\em(z_2)
+\frac{1-q^2}{1+q^2} \oint \dzzz{}\fp^2(z)\ot\em^2(z)}
$$
which obviously coincide with \r{int2}  for $n=2$.

As we can see from this exercise the most subtle point of these procedure is
the vanishing of the cross terms \r{cros}. We can observe these vanishing
property in general case. Indeed, the relations
$$
\kr_{\pm,\mp}= (\Pfpmd\ot 1)\ov\kr
=(1\ot \Pemp)\ov\kr=(\Pfpmd\ot \Pemp)\ov\kr
$$
are equivalent to the equalities
$$
\oint\cdots \oint \prod\dzz{i}\sk{f(z_1)\cdots f(z_n)}_\pm\ot
\sk{e(z_1)\cdots e(z_n)}_\mp=$$
\beq \label{ur1}
=\oint\cdots \oint \prod\dzz{i}\sk{f(z_1)\cdots f(z_n)}_\pm
\ot e(z_1)\cdots e(z_n)=\eeq
\beq\label{ur2}
=\oint\cdots \oint \prod\dzz{i}f(z_1)\cdots f(z_n)
\ot \sk{e(z_1)\cdots e(z_n)}_\mp\ .
\eeq

To proceed further we need the following Proposition \ref{Khf}.
Let $I\subset \{1,,...,n\}$, $I=\{ i_1<...<i_k\}$,
 $J=\{1,...,n\}\setminus I$, $J=\{ j_1<...<j_{n-k}\}$ be
some subsets from the set $\{1,...,n\}$.
For these sets
denote by $g_{I,J}(z_1,...,z_n)$ and $g'_{I,J}(z_1,...,z_n)$
the following formal power series:
$$
 g_{I,J}(z_1,...,z_n)=\prod_{i\in I,j\in J\atop i>j}
g\left(\frac{z_i}{z_j}\right),
\qquad
 g'_{I,J}(z_1,...,z_n)=\prod_{i\in I,j\in J\atop i>j}
g'\left(\frac{z_j}{z_i}\right).$$

\begin{prop}\label{Khf}
The following equalities of formal power series take place in
$\ov U_f^>\ot\ov U_e^<$:
\bse\label{Kf1}
\be\label{Kf1a}
f(z_1)\cdots f(z_n)=
\sum_{I,\ J}
g_{I,J}(z_1,...,z_n)
\left(f(z_{i_1})\cdots f(z_{i_k})\right)_+
\left(f(z_{j_1})\cdots f(z_{j_{n-k}})\right)_-\ ,
\ee
\be\label{Kf1b}
e(z_1)\cdots e(z_n)=\sum_{I,\ J}
g'_{I,J}(z_1,...,z_n)
\left(e(z_{i_1})\cdots e(z_{i_k})\right)_-
\left(e_(z_{j_1})\cdots e(z_{j_{n-k}})\right)_+\ .
\ee
\ese
\end{prop}
{\it Proof.}
For the proof let us apply \r{addition} to the product $f(z_1)\cdots
f(z_n)$.
We obtain
$$f(z_1)\cdots f(z_n)=\mu
(\Pfpd\ot \Pfm)(1\ot f(z_1)+f(z_1)\ot \psi^+(z_1))\cdots
(1\ot f(z_n)+f(z_n)\ot \psi^+(z_n))=
$$
$$
=\sum_{I\subset \{1,...,n\}, I=\{i_1<...<i_k\},}
\left(f(z_{i_1})\cdots f(z_{i_k})\right)_+
\left(a(z_1)\cdots a(z_n)\right)_-,$$
where $a(z_k)=\psi^+(z_k)$ if $k\in I$ and
$a(z_k)=f(z_k)$ otherwise. Moving all the $\psi^+(z_k)$ to the left
in the product $a(z_1)\cdots a(z_n)$
 and noting that $\varepsilon(\psi^+(z_k))=1$, we get the statement
 of the Proposition.

Let us consider the sum of the integrals
$$\sum_{I,\ J}
\oint\cdots \oint \prod\dzz{i}
g_{I,J}(z_1,\ldots,z_n)
\sk{f(z_{i_1})\cdots f(z_{i_k})}_+\sk{f(z_{j_1})\cdots f(z_{j_{n-k}})}_-
\ot e(z_{1})\cdots e(z_{n})$$
which means that we replaced the product $f(z_1)\cdots f(z_n)$ in the first
tensor space using the statement of the proposition \ref{Khf}.

We can reorder the product $e(z_1)\cdots e(z_n)$ by means of \r{L31}
and rewrite this integral in a form
$$\sum_{I,\ J}\oint\cdots \oint \prod\dzz{i}
\sk{f(z_{i_1})\cdots f(z_{i_k})}_+\sk{f(z_{j_1})\cdots f(z_{j_{n-k}})}_-
\ot
e(z_{i_1})\cdots e(z_{i_k})e(z_{j_1})\cdots e(z_{j_{n-k}})$$
Now we separate integrations over $z_{i_1},...,z_{i_k}$ and over
$z_{j_1},...,z_{j_{n-k}}$ and for the first integration use \r{ur1} $(+)$
and for the second -- \r{ur1} $(-)$. As a result, we have
$$\sum_{I,\ J}
\oint\cdots \oint \prod\dzz{i}g_{I,J}(z_1,\ldots,z_n)
\sk{f(z_{i_1})\cdots f(z_{i_k})}_+\sk{f(z_{j_1})\cdots f(z_{j_{n-k}})}_-
\ot
e(z_{1})\cdots e(z_{n})=$$
$$=\sum_{I,\ J}\oint\cdots \oint \prod\dzz{i}
\sk{f(z_{i_1})\cdots f(z_{i_k})}_+\sk{f(z_{j_1})\cdots f(z_{j_{n-k}})}_-
\ot
\sk{e(z_{i_1})\cdots e(z_{i_k})}_-\sk{e(z_{j_1})\cdots e(z_{j_{n-k}})}_+.$$
 This proves that the cross terms in the integral
$$
\oint\cdots \oint \prod\dzz{i} \sk{f(z_1)\cdots f(z_{n})}_+
\ot \sk{e(z_{1})\cdots e(z_{n})}_-
$$
 vanish.

\section*{Acknowledgment}

\qquad
Authors are indebted to A.~Okounkov for the proof of the identity
\r{l-stat} and to S. Kharchev
for useful discussions on the first stages of this work.

S.~Pakuliak also thanks the Physikalisches Institut,
Universit\"at Bonn;
CERN Theory Division and  Universit\'e d'Angers for hospitality.

The work of S.Kh. and of S.P.
 was supported in part by the grants INTAS OPEN 97-01312,
RFBR-CNRS grant PICS N 608/RFBR 98-01-22033.
S.Kh. was supported also  by RFBR grant 98-01-00303 and S.P.
by grants RFBR 97-01--1041  and grant of Heisenberg-Landau
program HL-99-12.

\app{Pairing calculations}
\label{normal}

The aim of this Appendix is two-fold. First, we will prove the
equivalence between equalities \r{con3} and combinatorial identity
\r{l-stat}. Second, we  explain several hints how to work with
orthogonal dual bases in nilpotent subalgebras of $\Uqdva$ during
this proof. At the end, we  exploit the evaluation homomorphism
in order to verify the factorization of the element $\ov\kr$.

\noindent{\bf 1. Equivalence between \r{con3} and \r{l-stat}.}

Using multiplicative presentation of the element $\ov\kr$ \r{du10}
we can define the dual bases in nilpotent  subalgebras
$\ov U_f^>$ and $\ov U_e^<$
corresponding to decomposition of the algebra $\Uqdva$ into two Borel
subalgebras associated with first Drinfeld Hopf structure:
\be\label{n13}
E^{\{p\}}=\left\{\left.
e_{p_1}e_{p_2}\cdots e_{p_n}
\ \right|\
{p_1\leq p_2\leq \cdots\leq p_n}
\right\},\\
F_{\{p\}}=\left\{\left.
C_{(\{p\})}(q)f_{-p_1}f_{-p_2}\cdots f_{-p_n}
\ \right| \ {p_1\leq p_2\leq \cdots\leq p_n}
\right\}\ ,
\\
\<F_{\{p\}},E^{\{m\}}\>=\delta_{\{p\}}^{\{m\}}=
\delta^{m_1}_{p_1}\delta^{m_2}_{p_2}\cdots\delta^{m_n}_{p_n}
,
\quad p_k,m_k\in\ZZ\ ,
\ee
where $C_{\{p\}}$ are given by \r{combinat1}.
We see from \r{n13}
that ordered monomials of the generators $e_n$ in the completion
$\ov{U}^<_e$ are dual to
the non-ordered monomials constructed from modes $f_n$ in the completion
$\ov{U}^>_f$.

Let us calculate left hand side of identity \r{con3a}.
To do this we rewrite the products of the currents
$f(z_1)\cdots f(z_k)$ and $e(z_1)\cdots e(z_k)$ in ordered form
using the dual basis
in the subalgebras $\ov{U}^>_f$ and  $\ov{U}^<_e$:
\bse\label{dp3}
\be\label{dp3a}
f(z_1)\cdots f(z_k)=\sum_{p_1\geq\ldots\geq p_k}
C_{\{p\}}(q)\<f(z_1)\cdots f(z_k),e_{-p_1}\cdots e_{-p_k}\>
f_{p_1}\cdots f_{p_k}\ ,
\ee
\be\label{dp3b}
e(z_1)\cdots e(z_k)=\sum_{s_1\leq\ldots\leq s_k}
C_{\{s\}}(q)\<f_{-s_1}\cdots f_{-s_k},e(z_1)\cdots e(z_k)\>
f_{p_1}\cdots f_{p_k}\ .
\ee
\ese
Then we have
\be\label{cal10}
\oint\cdots\oint \dzz{1}\cdots\dzz{k}\
\sk{\skk{f(z_1)\cdot f(z_k)}_- \fp(z) }_+\ot
\sk{\skk{e(z_1)\cdot e(z_k)}_+ \em(z) }_-=\\
=
\sk{\sk{
\sum_{p_1\geq\ldots\geq p_k}C_{\{p\}}(q)\ f_{p_1}\cdots f_{p_k}
}_-\fp(z)}_+\ot
\sk{\sk{
\sum_{s_1\leq\ldots\leq s_k}C_{\{s\}}(q)\ e_{s_1}\cdots e_{s_k}
}_+\em(z)}_-\times\\
\times
\oint\cdots\oint \dzz{1}\cdots\dzz{k}\
\<f_{-s_1}\cdots f_{-s_k},e(z_1)\cdots e(z_k)\>
\<f(z_1)\cdots f(z_k), e_{-p_1}\cdots e_{-p_k}\>=\\
\rav{m22}
\sk{
\sum_{0\geq p_1\geq\ldots\geq p_k}C_{\{p\}}(q)\ f_{p_1}\cdots f_{p_k}
\fp(z)}_+\ot
\sk{
\sum_{0\leq s_1\leq\ldots\leq s_k}C_{\{s\}}(q)\ e_{s_1}\cdots e_{s_k}
\em(z)}_-\times\\
\times
k!(q^{-1}-q)^k
\<f_{-s_1}\cdots f_{-s_k}, e_{-p_1}\cdots e_{-p_k}\>=\\
\rav{m3}
k!
\sum_{0\leq p_1\leq\ldots\leq p_k}C_{\{p\}}(q)\
\sk{ f_{-p_1}\cdots f_{-p_k}
\fp(z)}_+\ot
\sk{e_{p_1}\cdots e_{p_k} \em(z)}_-\ .
\ee

The last step to prove the equivalence of \r{con3a} to the combinatorial formula
\r{l-stat} is to use
the formulas
\be\label{auxid}
\sk{ f_{-p_1}\cdots f_{-p_k}
\fp(z)}_+=\prod_{m=1}^k(1-q^{2m})q^{2(p_k+2p_{k-1}+\cdots+kp_1)}
\fp^{k+1}(z)z^{-p_1-\cdots-p_k}\ ,\\
\sk{ e_{p_1}\cdots e_{p_k}
\em(z)}_-=\prod_{m=1}^k(1-q^{2m})q^{2(p_k+2p_{k-1}+\cdots+kp_1)}
\em^{k+1}(z)z^{p_1+\cdots+p_k}\
\ee
which are
 particular cases of more general formulas
\bse\label{prov}
\be\label{prova}
\sk{\fm^k(z_1)\fp^{n-k+1}(z_2)}_+
=\prod_{m=n-k-1}^{n}
\frac{1-q^{2m}}{1-q^{2m}z_1/z_2}\fp^{n+1}(z_2)\ ,
\ee
\be\label{provc}
\sk{\ep^k(z_1)\em^{n-k+1}(z_2)}_-=\prod_{m=n-k-1}^{n}
\frac{1-q^{2m}}{1-q^{2m}z_2/z_1}\em^{n+1}(z)\ ,
\ee
\be\label{provb}
\sk{\fm^{n-k+1}(z_1)
\fp^k(z_2)}_-=\prod_{m=n-k-1}^{n}
\frac{(1-q^{2m})z_1/z_2}{1-q^{2m}z_1/z_2}\fm^{n+1}(z_1)\ ,
\ee
\be\label{provd}
\sk{\ep^{n-k+1}(z_1)\em^k(z_2)}_+=
\prod_{m=n-k-1}^{n}\frac{(1-q^{2m})z_2/z_1}{1-q^{2m}z_2/z_1}
\ep^{n+1}(z_1)\ .
\ee
\ese

\noindent{\bf 2. Equivalence between multiplicative and integral
forms of the pairing tensor.}

In $\overline{U}^>_f\ot \overline{U}^<_e$  there is an identity
\bn
\overrightarrow{\prod_{n\in\ZZ}}
\exp_{q^2}\left((q^{-1}-q)f_{-n}\ot e_{n}\right)=
:\exp\sk{(q^{-1}-q)\oint \underline{dz}\ f(z)\otimes e(z)}:\ ,
\label{en}
\ed
where, as usual, the dots in r.h.s. means the antiordering in left tensor
space and ordering in right tensor space.

In order to prove \r{en} we need
two formulas
\be\label{m22}
\oint\cdots\oint\ \underline{dz_n}\ldots \underline{dz_n}
\ \langle f(w_1)\cdots f(w_n), e(z_1)\cdots e(z_n) \rangle\
\langle f(z_1)\cdots f(z_n), e(u_1)\cdots e(u_n) \rangle
=\\
= n! (q^{-1}-q)^{-n}
\langle f(w_1)\cdots f(w_n), e(u_1)\cdots e(u_n) \rangle
\ee
and
\be\label{m3}
\langle f_{-p_1}f_{-p_2}\cdots
f_{-p_n}, e_{s_1}e_{s_2}\cdots e_{s_n} \rangle
= \delta_{p_1,s_1}\delta_{p_2,s_2}
\cdots \delta_{p_n,s_n} C^{-1}_{\{p\}}(q) (q^{-1}-q)^{-n}
\ee
for $p_1\leq p_2\leq\ldots\leq p_n$ and $s_1\leq s_2\leq\ldots\leq s_n$.
The first formula can be proved using simple combinatorics and explicit
formulas for the pairing \r{pr} and the second one is the direct
consequence
of the factorization property of this pairing.

Now the statement of \r{en}
follows from the following calculation
\be\label{m4}
\frac{1}{n!}
\oint\cdots\oint\ \underline{dz_1}\ldots \underline{dz_n}
\ f(z_1)\cdots f(z_n)\ot e(z_1)\cdots e(z_n) =\\
\rav{dp3}
\frac{1}{n!}\sum_{p_1\leq \cdots\leq p_n\atop s_1\leq \cdots\leq s_n}
C_{\{p\}}(q) C_{\{s\}}(q)
f_{-s_1}\cdots f_{-s_n}\ot e_{p_1}\cdots e_{p_n} \times\\
\times
\oint\cdots\oint\ \underline{dz_n}\ldots \underline{dz_n}
\ \langle f_{-p_1}\cdots f_{-p_n}\ot e(z_1)\cdots e(z_n) \rangle\
\langle f(z_1)\cdots f(z_n)\ot e_{s_1}\cdots e_{s_n} \rangle =\\
\ravv{m22}{m3}
\sum_{p_1\leq p_2\leq \ldots \leq p_n} C_{\{p\}}(q)
f_{-p_1}\cdots f_{-p_n}\ot e_{p_1}\cdots e_{p_n}
\ee

\noindent{\bf 3. Evaluation map and universal $\kr$-matrix.}

Let $E$, $F$ and $H$ be generators of $U_q(\frak{sl}_2)$ with standard
commutation relations
\be\label{ev-eq}
q^H Eq^{-H}=q^2E,\quad q^H Fq^{-H}=q^{-2}F,\quad
[E,F]=\frac{q^H-q^{-H}}{q-q^{-1}}
\ee
An evaluation map $\Ev_a:\Uqdva\to U_q(\frak{sl}_2)\otimes
\CC[[a,a^{-1}]]$ is defined as follows:
\be\label{ev-def}
\Ev_a(e_n)=a^nq^{nH}E,\quad \Ev_a(f_n)=a^nFq^{nH},\quad
\Ev_a(\psi_n^\pm)=\pm(q-q^{-1})a^nq^{nH}\sk{EF-q^{2n}FE}.
\ee
Under this map the half currents will take the form
\be\label{eval}
\Ev_a\sk{f_-(z)}=-F\frac{1}{1-\frac{z}{a}q^{-H}}\ ,\qquad
\Ev_a\sk{e_+(z)}=\frac{1}{1-\frac{a}{z}q^{H}}\ E\ ,\\
\Ev_a\sk{f_+(z)}=F\frac{\frac{a}{z}q^{H}}{1-\frac{a}{z}q^{H}}\ ,\qquad
\Ev_a\sk{e_-(z)}=-\frac{\frac{z}{a}q^{H}}{1-\frac{z}{a}q^{-H}}\ E\ .
\ee
Using formulas \r{powers} we can write down the formulas for the projections
of two-tensors $I^{(n)}_{\pm,\mp}$ as the contour integrals over
unit circles ($|a|<1$, $|a|<1$)
\bse\label{eval1}
\be\label{eval1a}
\Ev_a\ot\Ev_b\sk{I^{(n)}_{+,-}}=\frac{(q^{-1}-q)^{2n-1}}
{(2\pi i)^n[n]_q}
\oint_{|z|=1}\dzz{}\ F\frac{\frac{a}{z}q^{H}}{1-\frac{a}{z}q^{H}}
\cdots F\frac{\frac{a}{z}q^{H}}{1-\frac{a}{z}q^{H}}\otimes
\frac{\frac{z}{b}q^{H}}{1-\frac{z}{b}q^{-H}}\ E\cdots
\frac{\frac{z}{b}q^{H}}{1-\frac{z}{b}q^{-H}}\ E\ ,
\ee
\be\label{eval1b}
\Ev_a\ot\Ev_b\sk{I^{(n)}_{-,+}}=\frac{(q^{-1}-q)^{2n-1}}
{(2\pi i)^n[n]_q}
\oint_{|z|=1}\dzz{}\ F\frac{1}{1-\frac{a}{z}q^{-H}}
\cdots F\frac{1}{1-\frac{a}{z}q^{-H}}\otimes
\frac{1}{1-\frac{z}{b}q^{H}}\ E\cdots
\frac{1}{1-\frac{z}{b}q^{H}}\ E\ .
\ee
\ese
Commutation relations \r{ev-eq} and Cauchy theorem allows one
to calculate these integrals to obtain:
\bse\label{ev2}
\be\label{ev2a}
\Ev_a\ot\Ev_b\sk{I^{(n)}_{+,-}}=\frac{(q^{-1}-q)^{2n-1}}{[n]_q}
\sum_{j=1}^n
\prod_{i=1\atop i\neq j}^n\frac{q^{2(i-j)}}{1-q^{2(i-j)}}
\prod_{\ell=1}^n\frac{Zq^{2(\ell+j-1)}}{1-Zq^{2(\ell+j-1)}}
F^n\ot E^n\ ,
\ee
\be\label{ev2b}
\Ev_a\ot\Ev_b\sk{I^{(n)}_{-,+}}=\frac{(q^{-1}-q)^{2n-1}}{[n]_q}
F^n\ot E^n
\sum_{j=1}^n
\prod_{i=1\atop i\neq j}^n\frac{1}{1-q^{2(i-j)}}
\prod_{\ell=1}^n\frac{1}{1-Z^{-1}q^{2(\ell+j-1)}}\ ,
\ee
\ese
where
$$
Z=\frac{a}{b}q^{ H\ot 1 - 1\ot H }\ .
$$
For the other hand  the expressions for
$\Ev_a\ot\Ev_b\sk{\kr^{(n)}_{\pm,\mp}}$ can be found in \cite{KTS}:
\bse\label{ev3}
\be\label{ev3a}
\Ev_a\ot\Ev_b\sk{\kr^{(n)}_{+,-}}=\frac{(q^{-1}-q)^{n}}{(n)_{q^{-2}}!}
\sk{\prod_{i=1}^{n}\frac{q^{2i}Z}{1-q^{2i}Z}}
F^n\ot E^n\ ,
\ee
\be\label{ev3b}
\Ev_a\ot\Ev_b\sk{\kr^{(n)}_{-,+}}=\frac{(q^{-1}-q)^{n}}{(n)_{q^2}!}
F^n\ot E^n
\sk{\prod_{i=1}^{n}\frac{1}{1-q^{2i}Z^{-1}}} \ .
\ee
\ese
We checked for $n=2,3$ that the elements \r{ev3} satisfy recurrence
relations \r{r-cal} and \r{l-cal} with $I^{(n)}_{\pm,\mp}$ given by
\r{ev2}. Unfortunately, we did not managed to verify directly the
corresponding differential equations for
$\Ev_a\ot\Ev_b\sk{\kr_{\pm,\mp}(\tau)}$. It would be interesting to find
such a proof.

\end{document}